%% file: mehta.tex
\documentclass[12pt]{article}
\usepackage{makeidx}
\usepackage{latexsym}
\usepackage{amsfonts}
\usepackage{amsmath}
\usepackage{amsthm}
\usepackage{amssymb}
\usepackage{amscd}
\usepackage[dvips]{graphicx}
\usepackage{cases}
\usepackage{color}
\usepackage{eepic}
\usepackage{setspace}
\usepackage{calrsfs}
\usepackage{authblk}
\theoremstyle{definition}
\newtheorem{theorem}{Theorem}[section]
\newtheorem{prop}[theorem]{Proposition}
\newtheorem{lemma}[theorem]{Lemma}
\newtheorem{corollary}[theorem]{Corollary}

\newtheorem{remark}[theorem]{Remark}
\newtheorem{conjecture}[theorem]{Conjecture}

\newenvironment{demo}[1]{%
  \trivlist
  \item[\hskip\labelsep
        {\bf #1.}]
}{%
\hfill\qedsymbol
  \endtrivlist
}
 \makeatletter
    
    \@addtoreset{equation}{section}
  \makeatother
\renewcommand{\mathcal}{\mathrsfs}

\title{
A generalization of the Mehta-Wang determinant
and Askey-Wilson polynomials
}
\author[1]{Masao ISHIKAWA\thanks{
Research partially supported by CNRS, Institut Camille Jordan, UMR 5208.}}
\affil[1]{\small Department of Mathematics, Faculty of Education, University of the Ryukyus, Nishihara, Okinawa 901-0213, Japan,
{\tt ishikawa@edu.u-ryukyu.ac.jp}
}
\author[2]{Hiroyuki TAGAWA\thanks{
Partially supported by Grant-in-Aid for Scientific Research~(C)~23540017.}
}
\affil[2]{Department of Mathematics, Faculty of Education, Wakayama University, Sakaedani, Wakayama 640-8510, Japan,
{\tt tagawa@math.edu.wakayama-u.ac.jp}
}
\author[3]{Jiang ZENG}
\affil[3]{Institut Camille Jordan, Universit\'e Claude Bernard Lyon 1, 69622 Villeurbanne cedex, France,
{\tt zeng@math.univ-lyon1.fr}
}
\date{
}
\def\defterm#1{{\sl #1}\/}
\def\newterm#1{{\sl #1}\/}

\renewcommand\tilde{\widetilde}
\newcommand\Pf{\operatorname{Pf}}

\def\Res{\operatorname{Res}}

\def\qbinom#1#2{\left[{{#1} \atop {#2}}\right]_{q}}
\newcommand\rdots{\mathinner{\mkern1mu\raise0pt\vbox{\kern7pt\hbox{.}}
     \mkern2mu\raise4pt\hbox{.}\mkern2mu\raise8pt\hbox{.}\mkern1mu}}
\def\V#1#2#3#4{V^{{#1},{#2}}\left({#3};{#4}\right)}

\def\qbin#1#2{\left[#1 \atop #2\right]_{q}}
\def\adj{\operatorname{adj}}
\def\covered{\mathinner{\mkern1mu\raise0pt\vbox{\kern7pt\hbox{$<$}}
     \mkern-4mu\raise2pt\hbox{.}\mkern2mu}}
\def\covers{\mathinner{\mkern1mu\raise0pt\vbox{\kern7pt\hbox{$>$}}
     \mkern-12mu\raise2pt\hbox{.}\mkern8mu}}

\def\img{{\imath}} 
\def\Z{{\mathbb{Z}}}

\def\i{{\boldsymbol{i}}}
\def\j{{\boldsymbol{j}}}
\def\k{{\boldsymbol{k}}}
\def\x{{\boldsymbol{x}}}
\input roby4.tex

\begin{document}
%
%
%
\maketitle
\kern-15pt
\begin{abstract}
Motivated by the Gaussian symplectic ensemble, Mehta and Wang evaluated the $n\times n$ determinant $\det((a+j-i)\Gamma(b+j+i))$ in 2000.
 When $a=0$, Ciucu and Krattenthaler computed the associated Pfaffian $\Pf((j-i)\Gamma(b+j+i))$
with an application to the two dimensional dimer system in 2011.  Recently we have  generalized the latter Pfaffian formula  with a $q$-analogue
by replacing the Gamma function by the moment sequence of the little $q$-Jacobi polynomials.
On the other hand, Nishizawa has found  a $q$-analogue of the Mehta--Wang formula. Our  purpose is
to generalize both the Mehta-Wang and Nishizawa formulae
 by using the moment sequence of the little $q$-Jacobi polynomials.
It turns out that the corresponding determinant  can be evaluated explicitly in terms of the Askey-Wilson polynomials.
\end{abstract}
%
%
%
\input mehta01.tex
%
%
\input mehta02.tex
%
%
\input mehta03.tex
%
%
\input mehta04.tex
%
%
%
%
%
%
%
%
%
%
%
%
%
%
%
%
\par\medskip\noindent
{\large\bf Acknowledgments}
We are grateful to Professor Christian Krattenthaler for his encouragement at the initial stage of this project.
%
%
%
%
%
%

%
%
%
%
%
\end{document}

%% file: roby4.tex
\newdimen\Squaresize \Squaresize=20pt
\newdimen\thickness \thickness=1pt         
                                                    
\def\Square#1{\hbox{\vrule width \thickness
   \vbox to \Squaresize{\hrule height \thickness\vss                            
      \hbox to \Squaresize{\hss#1\hss}
   \vss\hrule height\thickness} 
\unskip\vrule width \thickness} 
\kern-\thickness}                                                            
                               
\def\vsquare#1{\vbox{\Square{$#1$}}\kern-\thickness}
\def\blank{\omit\hskip\Squaresize}

\def\fibyoung#1{\let\\=\cr              
\vbox{\smallskip\offinterlineskip
\halign{&\vsquare{##}\cr #1}}\,}

\def\borderlessrect#1#2{\hbox{
   \vbox to \Squaresize{\vskip \thickness \vss
      \hbox to #2 {\hss #1\hss}
   \vss\vskip\thickness} 
\unskip \hskip -\thickness
} 
\kern-\thickness}                                                            
                               
\def\vborderlessrect#1#2{\vbox{\borderlessrect{$#1$}{#2}}\kern-\thickness}

\def\borderless#1{\omit\vborderlessrect{#1}{\Squaresize}}
\def\borderlessrc#1#2{\omit\vborderlessrect{#1}{#2}}

\def\msquare#1{\vbox{\hbox{\vrule width \thickness
   \vbox to \Squaresize{\hrule height \thickness
      \hbox to \Squaresize{\hfil{\sevenrm #1}}
   \vfil\hrule height\thickness}
\unskip\vrule width \thickness}
\kern-\thickness}\kern-\thickness}

\def\twosquare#1#2{\vbox{\hbox{\vrule width \thickness 
   \vbox to \Squaresize{\hrule height \thickness
      \hbox to \Squaresize{\hfil{\rm #1}}\vss
      \hbox to \Squaresize{\hss{#2}\hss}
   \vfil\hrule height\thickness}
\unskip\vrule width \thickness}
\kern-\thickness}\kern-\thickness}

\def\twoblank#1#2{\vbox{\hbox{
   \vbox to \Squaresize{\vskip 2pt
      \hbox to \Squaresize{\hfil{\sevenrm #1}\ }\vss
      \hbox to \Squaresize{\hss{#2}\hss}
   \vfil}\unskip\kern-\thickness}
}\unskip\kern-\thickness}

\def\young#1{
\def\>{\blank}
\def\<{\borderless}
\def\*{\borderlessrc}
\def\p{\omit\msquare}
\def\t{\omit\twosquare}
\def\b{\omit\twoblank}
\let\\=\cr 
\vbox{\smallskip\offinterlineskip
\halign{&\vsquare{##}\cr #1}}}

\newdimen\smsquaresize \smsquaresize=12pt
\newdimen\smthickness \smthickness=.5pt
\font\smcellfont=cmss8 scaled \magstep0

\def\smsquare#1{\hbox{\vrule width \smthickness
   \unskip\vbox to \smsquaresize{\hrule height \smthickness\vss
      \hbox to \smsquaresize{\hss{\smcellfont #1}\hss}
   \vss\hrule height\smthickness} 
\unskip\vrule width \smthickness} 
\kern-\smthickness}

\def\smvsquare#1{\vbox{\smsquare{$#1$}}\kern-\smthickness}
\def\blank{\omit\hskip\smsquaresize}

\def\smyoung#1{\let\\=\cr 
\vbox{\smallskip\offinterlineskip
\halign{&\smvsquare{##}\cr #1}}}
\newdimen\vsmsquaresize \vsmsquaresize=10pt
\newdimen\vsmthickness \vsmthickness=.5pt
\font\vsmcellfont=cmsl8 scaled \magstep0
\font\vsmletterfont=cmr6 scaled \magstep0

\def\vsmsquare#1{\hbox{\vrule width \vsmthickness
   \unskip\vbox to \vsmsquaresize{\hrule height \vsmthickness\vss
      \hbox to \vsmsquaresize{\hss{\vsmcellfont #1}\hss}
   \vss\hrule height\vsmthickness} 
\unskip\vrule width \vsmthickness} 
\kern-\vsmthickness}
\def\vsmvsquare#1{\vbox{\vsmsquare{#1}}\kern-\vsmthickness}
\def\vsmblank{\omit\hskip\vsmsquaresize}
\def\vsmborderless#1{\hbox{\hskip \vsmthickness\unskip
   \vbox to \vsmsquaresize{\vss
      \hbox to \vsmsquaresize{\hss{\vsmletterfont #1}\hss}
   \vss} 
\unskip\hskip \vsmthickness} 
\kern-\vsmthickness}                                                            \def\vsmvborderless#1{\vbox{\vsmborderless{#1}}\kern-\vsmthickness}

\def\vsmyoung#1{
\def\>{\vsmblank}
\def\<{\omit\vsmvborderless}
\let\\=\cr 
\vbox{\smallskip\offinterlineskip
\halign{&\vsmvsquare{##}\cr #1}}}

\newdimen\edgesize \edgesize=20pt
\newdimen\eedgesize \eedgesize=21pt
\newdimen\doublesize \doublesize=41pt
\newdimen\ddoublesize \ddoublesize=43pt

\newdimen\triplesize \triplesize=62pt
\newdimen\tetrasize \tetrasize=85pt
\newdimen\futosa \futosa=1pt         

\def\Hako#1{\hbox{\vrule width \futosa
   \vbox to \eedgesize{\hrule height \futosa\vss                            
      \hbox to \edgesize{\hss#1\hss}
   \vss\hrule height\futosa} 
\unskip\vrule width \futosa} 
\kern-\futosa}                                                            
                               
\def\Seihokei#1{\vbox{\Hako{#1}}\kern-\futosa}

\def\Horizontal#1{\hbox{\vrule width \futosa
   \vbox to \eedgesize{\hrule height \futosa\vss                            
      \hbox to \doublesize{\hss#1\hss}
   \vss\hrule height \futosa} 
\unskip\vrule width \futosa} 
\kern-\futosa}                                                            
                               
\def\Vertical#1{\hbox{\vrule width \futosa
   \vbox to \ddoublesize{\hrule height \futosa\vss                            
      \hbox to \edgesize{\hss#1\hss}
   \vss\hrule height \futosa} 
\unskip\vrule width \futosa} 
\kern-\futosa}                                                            
                               
\def\NS#1{\vbox{\Hako{$#1$}\vskip22pt}\kern-\futosa}
\def\NSS#1{\vbox{\Hako{$#1$}\vskip42pt}\kern-\futosa}
\def\NSSS#1{\vbox{\Hako{$#1$}\vskip64pt}\kern-\futosa}
\def\H#1{\vbox{\Horizontal{$#1$}}\kern-\futosa}
\def\HH#1#2{\vbox{\Horizontal{$#1$}\Horizontal{$#2$}}\kern-\futosa}
\def\VER#1{\vbox{\Vertical{$#1$}}\kern-\futosa}
\def\VS#1{\vbox{\Vertical{$#1$}\vskip20pt}\kern-\futosa}
\def\VSS#1{\vbox{\Vertical{$#1$}\vskip42pt}\kern-\futosa}
\def\VV#1#2{\vbox{\Vertical{$#1$}\Vertical{$#2$}}\kern-\futosa}
\def\NV#1#2{\vbox{\Hako{$#1$}\Vertical{$#2$}}\kern-\futosa}
\def\NVS#1#2{\vbox{\Hako{$#1$}\Vertical{$#2$}\vskip22pt}\kern-\futosa}
\def\YTT#1#2#3{\vbox{\Horizontal{$#1$}\hbox{\vbox{\Vertical{$#2$}}\kern-\futosa\vbox{\Vertical{$#3$}}\kern-\futosa}}\kern-\futosa}

\def\domino#1{
\def\ns{\omit\NS}
\def\nss{\omit\NSS}
\def\nsss{\omit\NSSS}
\def\h{\omit\H}
\def\hh{\omit\HH}
\def\V{\omit\VER}
\def\Vs{\omit\VS}
\def\Vss{\omit\VSS}
\def\vsss{\omit\VSSS}
\def\vv{\omit\VV}
\def\nv{\omit\NV}
\def\nvs{\omit\NVS}
\def\ytt{\omit\YTT}
\let\\=\cr 
\vbox{\smallskip\offinterlineskip
\halign{&\Seihokei{##}\cr #1}}}

%% file: mehta01.tex
%
%
%
\section{Introduction and the main results}
%
%
%
%
%
%
Motivated by the Gaussian symplectic ensemble, 
Mehta and Wang \cite{MW} obtain the determinant identity
\begin{equation}
\det\bigl((a+j-i)\Gamma(b+i+j)\bigr)_{0\leq i,j\leq n-1}
= D_{n}\prod_{i=0}^{n-1}i!\Gamma(b+i),
\label{eq:Mehta-Wang}
\end{equation}
where (N.B. the binomial coefficient ${n\choose k}$ is missing in \cite[(7)]{MW})
\begin{align}
D_{n} 
=\sum_{k=0}^{n}
(-1)^{k}{n\choose k}\left(\frac{b-a}2\right)_{k}\left(\frac{a+b}2\right)_{n-k},
\end{align}
where $(\alpha)_n=\frac{\Gamma(\alpha+n)}{\Gamma(\alpha)}$ is known as the \defterm{rising factorial}.
This $D_{n}$ satisfies the three term recurrence relation
\begin{equation}
D_{-1}=0,\quad
D_{0}=1,\quad
D_{n+1}=aD_{n}+n(b+n-1)D_{n-1},
\label{eq:rec-D}
\end{equation}
which can be considered as the recurrence relation for a special case of 
the Meixner-Pollaczek polynomials (see \cite{MW,Ni}),
and one may notice that the sequence $\{\Gamma(b+n)\}_{n\geq0}$ of the Gamma functions in the left-hand side
can be considered as the moment sequence of the Laguerre polynomials
(see, for example, \cite{KLS,KS,Vie}).
Nishizawa \cite{Ni} obtains a $q$-analogue of \eqref{eq:Mehta-Wang},
which will be stated below.
Here we replace the Gamma functions by the moments of the little $q$-Jacobi
polynomials
and show that we obtain a special case of the Askey-Wilson polynomials as $D_n$,
which also generalize the two results in our previous papers \cite{ITZ1,ITZ2}.
Before we describe our results we need more notation.
\par\smallskip
%
%
%
%
%
%
Throughout this paper we use the standard  notation for $q$-series (see \cite{AAR,GR,KLS,KS}):
\begin{equation*}
    (a;q)_{\infty}=\prod_{k=0}^{\infty}(1-aq^{k}),\qquad
    (a;q)_{n}=\frac{(a;q)_{\infty}}{(aq^{n};q)_{\infty}}
\end{equation*}
for any integer $n$.
Usually  $(a;q)_{n}$ is called the  \defterm{$q$-shifted factorial},
and we frequently use the compact notation:
\begin{align*}
    &(a_{1},a_{2},\dots,a_{r};q)_{n}=(a_{1};q)_{n}(a_{2};q)_{n}\cdots(a_{r};q)_{n}.
\end{align*}
The \defterm{${}_{r+1}\phi_{r} $ basic hypergeometric series} is defined by
\begin{align}
{}_{r+1}\phi_{r}\left[\,
{{a_{1},a_{2},\dots,a_{r+1}}\atop{b_{1},\dots,b_{r}}};q,z
\,\right]
=\sum_{n=0}^{\infty}\frac{(a_{1},a_{2},\dots,a_{r+1};q)_{n}}{(q,b_{1},\dots,b_{r};q)_{n}}z^{n}.
\label{eq:phi-def}
\end{align}
Here we also use the $q$-Gamma function
\[
\Gamma_{q}(z)=(1-q)^{1-z}\frac{(q;q)_{\infty}}{(q^z;q)_{\infty}},
\]
the $q$-integer $[n]_{q}=\frac{1-q^n}{1-q}$ and the $q$-factorial $[n]_{q}!=\prod_{k=1}^n[k]_{q}$.
By taking the limit $q\rightarrow1$, we obtain
the hypergeometric function
\[
{}_{r+1}F_{r}\left[\,
{{\alpha_{1},\alpha_{2},\dots,\alpha_{r+1}}\atop{\beta_{1},\dots,\beta_{r}}};z
\,\right]
=\sum_{n=0}^{\infty}\frac{(\alpha_{1})_{n}(\alpha_{2})_{n}\cdots(\alpha_{r+1})_{n}}{n!\,(\beta_{1})_{n}\cdots(\beta_{r})_{n}}z^{n}.
\]
%
%
%
%
%
%
The \defterm{Askey-Wilson polynomials} $p_{n}(x)$ \cite{GR,KLS,KS} satisfy
the well-known recurrence relation
\begin{align}
2x p_{n}(x)=A_{n}p_{n+1}(x)+B_{n}p_{n}(x)+C_{n}p_{n-1}(x),
\quad n\geq0,
\label{eq:AW-rec}
\end{align}
with $p_{-1}(x)=0$, $p_{0}(x)=1$, where
\begin{align*}
A_{n}
&=\frac{1-abcdq^{n-1}}{\left(1-abcdq^{2n-1}\right)\left(1-abcdq^{2n}\right)},
\\
C_{n}
&=\frac{\left(1-q^{n}\right)\left(1-abq^{n-1}\right)\left(1-acq^{n-1}\right)\left(1-adq^{n-1}\right)}
{\left(1-abcdq^{2n-2}\right)\left(1-abcdq^{2n-1}\right)}
\\&\quad\times
\left(1-bcq^{n-1}\right)\left(1-bdq^{n-1}\right)\left(1-cdq^{n-1}\right),
\end{align*}
and
\begin{align*}
B_{n}&=a+a^{-1}
-A_{n}a^{-1}\left(1-abq^{n}\right)\left(1-acq^{n}\right)\left(1-adq^{n}\right)
\\&\quad
-C_{n}a/\left(1-abq^{n-1}\right)\left(1-acq^{n-1}\right)\left(1-adq^{n-1}\right).
\end{align*}
They have the basic hypergeometric expression
\begin{equation}
p_{n}(x;a,b,c,d;q)=\frac{(ab,ac,ad;q)_{n}}{a^n}\,
{}_{4}\phi_{3}\left(
{{q^{-n},abcdq^{n-1},ae^{\img\theta},ae^{-\img\theta}}\atop{ab,ac,ad}};q,q
\right)
\label{eq:Askey-Wilson}
\end{equation}
with $x=\cos\theta$, where $\img=\sqrt{-1}$.
%
%
%
We also use the symbol
\[
\chi(A)=\begin{cases}1&\text{ if $A$ is true,}\\0&\text{ if $A$ is false.}\end{cases}
\]
\par\smallskip
%
%
%
%
%
%
%
%
In \cite{ITZ1}
we have proven the Hankel determinant identity
\begin{align}
&\det\left(\frac{(aq;q)_{i+j+r-2}}{(abq^2;q)_{i+j+r-2}}\right)_{1\leq i,j\leq n}
  =a^{\frac{n(n-1)}2}
q^{\frac{n(n-1)(2n-1)}6+\frac{n(n-1)r}2}
  \nonumber\\&\qquad\qquad\times
  \prod_{k=1}^{n}\frac{(q,bq;q)_{k-1}(aq;q)_{k+r-1}}
  {(abq^{2};q)_{k+n+r-2}}
  \label{eq:q-general-hankel}
\end{align}
for a positive integer $n$.
Here
\[
\mu_{n}=\frac{(aq;q)_{n}}{(abq^2;q)_{n}}\quad(n=0,1,2,\dots)
\]
is the moments of the little $q$-Jacobi polynomials.
%
%
%
%
%
%
In our previous paper \cite{ITZ2} we have exploited the Pfaffian identity
\begin{align}
&\Pf\left((q^{i-1}-q^{j-1})\frac{(aq;q)_{i+j+r-2}}{(abq^2;q)_{i+j+r-2}}\right)_{1\leq i,j\leq 2n}
\nonumber\\&=
a^{n(n-1)}q^{\frac{n(n-1)(4n+1)}3+n(n-1)r}\prod_{k=1}^{n-1}(bq;q)_{2k}
\prod_{k=1}^n\frac{(q;q)_{2k-1}(aq;q)_{2k+r-1}}
{(abq^2;q)_{2(k+n)+r-3}}
\label{eq:pf-special}
\end{align}
for a positive integer $n$
(see also \cite{Las1,Las2}).
\par\smallskip
%
%
%
%
In \cite{Ni} Nishizawa has proven the $q$-analogue of the Mehta-Wang result:
\begin{align}
&\det\bigl([a+j-i]_{q}\Gamma_{q}(b+i+j)\bigr)_{0\leq i,j\leq n-1}
\nonumber\\
&\qquad
= q^{na+n(n-1)b/2+n(n-1)(2n-7)/6}\,D_{n,q}
\prod_{k=0}^{n-1}
[k]_{q}!\cdot\Gamma_{q}(b + k),
\label{eq:Nishizawa}
\end{align}
where $D_{n,q}$ satisfies the recurrence relation
\begin{align}
&D_{-1,q}=0,\quad D_{0,q}=1,\quad
D_{n+1,q}=q^{-a+n}[a]_{q}D_{n,q}+q^{-a-b}[n]_{q}[b+n-1]_{q}D_{n-1,q}.
\end{align}
%
%
Comparing this recurrence relation with the recurrence equation
\begin{equation}
2xQ_n(x)=Q_{n+1}(x)+(A+B)q^nQ_n(x)+(1-q^n)(1-ABq^{n-1})Q_{n-1}(x)
\label{eq:rec-AlSalam-Chihara}
\end{equation}
of the Al-Salam--Chihara polynomials
\begin{equation}
Q_{n}(x)=Q_{n}(x;A,B;q)
=\frac{(AB;q)_{n}}{A^n}
{}_{3}\phi_{2}\left({{q^{-n},Ae^{\img\theta},Ae^{-\img\theta}}\atop{AB,0}}\,;\,q;q\right)
\label{eq:def-AlSalam-Chihara}
\end{equation}
with $x=\cos\theta$ (see \cite{KLS,KS}),
we may remark that $D_{n,q}$ can be considered as a special case 
of the Al-Salam--Chihara polynomials because
\begin{equation}
D_{n,q}=(-\img)^{n}q^{-\frac{a+b}2n}(1-q)^{-n}Q_{n}\left(0;q^{\frac{a+b}2}\img,-q^{\frac{b-a}2}\img;q\right).
\label{eq:D-by-AlSalamChihara}
\end{equation}
By this observation,
we can write $D_{n,q}$ explicitly as
\begin{align}
D_{n,q}=
\frac{(q^b;q)_n}{q^{n(a+b)}(q-1)^{n}} \sum_{k=0}^nq^k\frac{(q^{-n};q)_k}{(q;q)_k}
 \prod_{j=0}^{k-1} \frac{1-q^{a+b+2j}}{1-q^{b+j}}.
\end{align}
%
%
%
One natural question we may ask is what can we obtain if we replace the $q$-Gamma function in the determinant of \eqref{eq:Nishizawa}
by the moment of the little $q$-Jacobi polynomials.
The aim of this paper is to answer this question,
and we can express the determinant by the Askey-Wilson polynomials.
%
%
%
\begin{theorem}
\label{th:corthm2}
Let $a$, $b$ and $c$ be parameters,
and let $n\geq1$ and $r$ be integers.
Then we have
\begin{align}
&\det\left(
(q^{i-1}-cq^{j-1})\frac{(aq;q)_{i+j+r-2}}{(abq^2;q)_{i+j+r-2}}
\right)_{1\leq i,j\leq n}\nonumber\\
&=
(-1)^na^{\frac{n(n-3)}{2}}q^{\frac{n(n+1)(2n-5)}{6}+\frac{n(n-3)r}{2}}(abcq^{r+1};q^2)_n
\prod_{k=1}^{n}\frac{(q;q)_{k-1}(aq;q)_{k+r}(bq;q)_{k-2}}{(abq^2;q)_{k+n+r-2}}
\nonumber\\
&\quad\times
{}_4\phi_3
\left(
\begin{matrix}
q^{-n}, 
a^{\frac{1}{2}}c^{\frac{1}{2}}q^{\frac{r+1}{2}}, 
-a^{\frac{1}{2}}c^{\frac{1}{2}}q^{\frac{r+1}{2}}, 
abq^{n+r}\cr
aq^{r+1}, 
a^{\frac{1}{2}}b^{\frac{1}{2}}c^{\frac{1}{2}}q^{\frac{r+1}{2}}, 
-a^{\frac{1}{2}}b^{\frac{1}{2}}c^{\frac{1}{2}}q^{\frac{r+1}{2}}
\end{matrix}; q,q
\right)\label{eqtnabcn3}\\
&=
(-\img)^na^{\frac{n(n-2)}{2}}c^{\frac{n}{2}}q^{\frac{n(n-2)(2n+1)}{6}+\frac{n(n-2)r}{2}}
\prod_{k=1}^{n}\frac{(q;q)_{k-1}(aq;q)_{k+r-1}(bq;q)_{k-2}}{(abq^2;q)_{k+n+r-2}}
\nonumber\\
&\quad\times
p_n\left(0;a^{\frac{1}{2}}c^{\frac{1}{2}}q^{\frac{r+1}{2}}\img,
-a^{\frac{1}{2}}c^{-\frac{1}{2}}q^{\frac{r+1}{2}}\img,
b^{\frac{1}{2}}\img,
-b^{\frac{1}{2}}\img;q\right).
\label{eqtnabcn3aw}
\end{align}
\end{theorem}
%
%
%
%
\begin{remark}
If we put $c=0$ in \eqref{eqtnabcn3}, then we recover 
our previous result \eqref{eq:q-general-hankel} easily by using the $q$-Chu-Vandermonde formula \cite[(1.5.3)]{GR}
\begin{equation}
{}_{2}\phi_{1}\left(
{{a,q^{-n}}\atop{c}}\,;\,q,q\right)
=\frac{(c/a;q)_{n}}{(c;q)_{n}}a^{n}.
\label{eq:Chu-Vandermonde}
\end{equation}
If we put $a=q^{\alpha-1}$, $b=0$, $c=q^{\gamma}$ and $r=0$ in \eqref{eqtnabcn3}, 
then the left-hand side equals
\[
\frac{q^{\frac{n(n-1)}2}(1-q)^{n^2}}{\{\Gamma_{q}(\alpha)\}^{n}}
\det\left([\gamma+j-i]_{q}\Gamma_{q}(\alpha+i+j-2)\right)_{1\leq i,j\leq n}
\]
because of $(q^\alpha;q)_{n}=(1-q)^{n}\cdot\frac{\Gamma_{q}(\alpha+n)}{\Gamma_{q}(\alpha)}$,
and the right-hand side equals 
\begin{align*}
(-\img)^nq^{\frac{n(n-2)}2\alpha+\frac{n}2\gamma+\frac{n(n-1)(n-2)}3}\prod_{k=1}^{n}(q;q)_{k-1}(q^\alpha;q)_{k-1}
\cdot Q_{n}(0;q^{\frac{\alpha+\gamma}2}\img,-q^{\frac{\alpha-\gamma}2}\img;q)
\end{align*}
because of the relation $Q_{n}(x;A,B;q)=p_{n}(x;A,B,0,0;q)$
between the Al-Salam--Chihara polynomials and the Askey-Wilson polynomials.
Hence we obtain Nishizawa's formula \eqref{eq:Nishizawa} as a corollary.
\end{remark}
In Section~\ref{sec:proof-theorems} we derive the following corollary from Theorem~\ref{th:corthm2}.
%
%
\begin{corollary}
\label{th:main}
Let $a$, $b$ and $c$ be parameters,
and let $n\geq1$ and $r$ be integers.
\begin{enumerate}
\item[(i)]
If the size $n=2m$ of the matrix is even,
then we have
\begin{align}
&\det\left((q^{i-1}-cq^{j-1})\frac{(aq;q)_{i+j+r-2}}{(abq^2;q)_{i+j+r-2}}\right)_{1\leq i,j\leq 2m}
\nonumber\\
&=
a^{2m(m-1)}c^mq^{\frac{2m(m-1)(4m+1)}{3}+2m(m-1)r}
\prod_{k=1}^{m}\left\{\frac{(q;q)_{2k-1}(aq;q)_{2k+r-1}(bq;q)_{2k-2}}
{(abq^2;q)_{2(k+m)+r-3}}\right\}^{2}
\nonumber\\
&\quad\times
{}_4\phi_3
\left(
\begin{matrix}
q^{-2m},b^{-1}q^{-2m+1},c,c^{-1}\cr
q,aq^{r+1},a^{-1}b^{-1}q^{1-4m-r}\cr
\end{matrix};
q^2,q^2
\right)
\label{eq:main-even}
\\
&=(-1)^m a^{m(2m-1)}b^m c^m q^{\frac{m (8m^2+3\,{m}-2)}{3}+m(2m-1)r}
\prod_{k=1}^{2m}\frac{(q;q)_{k-1}(aq;q)_{k+r-1}}
{(abq^2;q)_{k+2m+r-2}}
\nonumber\\
&\quad\times
\prod_{k=1}^{m}
\{(bq;q)_{2k-2}\}^2\cdot
 p_{m}\left(\frac{c+c^{-1}}2;1,q,aq^{r+1},a^{-1}b^{-1}q^{1-4m-r};q^2\right).
\label{eq:main-even2}
\end{align}
\item[(ii)]
If the size $n=2m+1$ of the matrix is odd,
then we have
\begin{align}
&\det\left((q^{i-1}-cq^{j-1})\frac{(aq;q)_{i+j+r-2}}{(abq^2;q)_{i+j+r-2}}\right)_{1\leq i,j\leq 2m+1}
\nonumber\\
&=
a^{2m^2}c^{m}q^{\frac{2m(m+1)(4m-1)}{3}+2m^2r}\cdot\frac{1-c}{1-q}\cdot
\prod_{k=1}^{m+1}\frac{
(q;q)_{2k-1}
(aq;q)_{2k+r-2}
(bq;q)_{2k-2}
}
{
(abq^{2};q)_{2(k+m-1)+r}}
\nonumber\\
&\quad\times
\prod_{k=1}^{m}\frac{(q;q)_{2k-1}(aq;q)_{2k+r}(bq;q)_{2k-2}}{(abq^2;q)_{2(k+m-1)+r}}
\cdot
{}_4\phi_3
\left(
\begin{matrix}
q^{-2m},b^{-1}q^{-2m+1},cq,c^{-1}q\cr
q^3,aq^{r+2},a^{-1}b^{-1}q^{-4m-r}
\end{matrix};
q^2,q^2
\right)
\label{eq:main-odd}
\\&=
(-1)^m a^{m(2m+1)}b^{m} c^m (1-c)
q^{\frac{m(8m^2+15m+4)}3+m(2m+1)r}
\prod_{k=1}^{2m+1}\frac{(q;q)_{k-1}(aq;q)_{k+r-1}}{(abq^2;q)_{k+2m+r-1}}
\nonumber\\&\quad\times
\prod_{k=1}^{m+1}(bq;q)_{2k-2}\cdot\prod_{k=1}^{m}(bq;q)_{2k-2}\cdot
 p_{m}\left(\frac{c+c^{-1}}2;q,q^2,aq^{r+1},a^{-1}b^{-1}q^{-4m-r-1};q^2\right).
\label{eq:main-odd2}
\end{align}
\end{enumerate}
\end{corollary}
%
%
%
%
%
\begin{remark}
If we put $c=1$ in \eqref{eq:main-even} for the even case,
then it is clear that the ${}_4\phi_3$ sum reduces to $1$,
so that the determinant becomes the product which equals the square
of the Pfaffian \eqref{eq:pf-special} obtained in \cite{ITZ2}.
Meanwhile, it does not suffice to prove \eqref{eq:pf-special}
since it is not so trivial to take the square root of the determinant
and determine the sign
(see \cite{CK,ITZ2}).
If we put $c=1$ in \eqref{eq:main-odd} for the odd case,
then the factor $(1-c)$ reduces the right-hand side to $0$.
\end{remark}
If we put $a=q^{\alpha}$, $b=q^{\beta}$ and $c=q^{\gamma}$
and let $q\rightarrow1$ in Theorem~\ref{th:corthm2},
then we obtain the following corollary.
%
%
%
\begin{corollary}
\label{cor:specilization}
Let $\alpha$, $\beta$ and $\gamma$ be parameters,
and let $n\geq1$ and $r$ be integers.
Then we have
\begin{align}
&\det\left((\gamma+j-i)\frac{(\alpha+1)_{i+j+r-2}}{(\alpha+\beta+2)_{i+j+r-2}}\right)_{1\leq i,j\leq n}
\nonumber\\&
=(-2)^{n}\left(\frac{\alpha+\beta+\gamma+r+1}2\right)_{n}\cdot
\prod_{k=1}^{n}\frac{(k-1)!(\alpha+1)_{k+r}(\beta+1)_{k-2}}
{(\alpha+\beta+2)_{k+n+r-2}}
\nonumber\\&\qquad\times
{}_3F_2
\left(
\begin{matrix}
-n,
\frac{\alpha+\gamma+r+1}2,
\alpha+\beta+n+r
\cr
\frac{\alpha+\beta+\gamma+r+1}2,
\alpha+r+1
\end{matrix}; 1
\right)
\label{eq:cor-special}
\\&
=(2\img)^{n}
\prod_{k=1}^{n}\frac{k!(\alpha+1)_{k+r-1}(\beta+1)_{k-2}}
{(\alpha+\beta+2)_{k+n+r-2}}
\nonumber\\&\qquad\times
\tilde{P}_{n}\left(0;\frac{\alpha+\gamma+r+1}2,\frac{\beta}2,\frac{\alpha-\gamma+r+1}2,\frac{\beta}2\right),
\label{eq:cor-special2}
\end{align}
where
\begin{align*}
&\tilde{P}_{n}(x;a,b,c,d)
\\&\quad
=\img^n\frac{(a+c)_{n}(a+d)_{n}}{n!}\cdot
{}_{3}F_{2}\left({{-n,n+a+b+c+d-1,a+\img x}\atop{a+c,a+d}};1\right)
\end{align*}
are the continuous Hahn polynomials (see \cite{KLS,KS}).
\end{corollary}
%
%
%
%
The \newterm{Wilson polynomials} $W_{n}(x;\alpha,\beta,\gamma,\delta)$ \cite{KLS,KS} are defined by
\begin{equation}
\frac{W_{n}(x^2;\alpha,\beta,\gamma,\delta)}{(\alpha+\beta)_{n}(\alpha+\gamma)_{n}(\alpha+\delta)_{n}}
={}_{4}F_{3}\left(
{{-n,\alpha+\beta+\gamma+\delta+n-1,\alpha+\img x,\alpha-\img x}\atop{\alpha+\beta,\alpha+\gamma,\alpha+\delta}};1
\right).
\end{equation}
%
%
%
%
%
%
By the same specialization as above,
we obtain the following corollary from Corollary~\ref{th:main}.
%
%
%
\begin{corollary}
\label{cor:wilson}
Let $\alpha$, $\beta$ and $\gamma$ be parameters,
and let $n\geq1$ and $r$ be integers.
\begin{enumerate}
\item[(i)]
If the size $n=2m$ of the matrix is even,
then we have
\begin{align}
&\det\left((\gamma+j-i)\frac{(\alpha+1)_{i+j+r-2}}{(\alpha+\beta+2)_{i+j+r-2}}\right)_{1\leq i,j\leq 2m}\nonumber\\
&=
\prod_{k=1}^{m}\left\{\frac{(2k-1)!(\alpha+1)_{2k+r-1}(\beta+1)_{2k-2}}
{(\alpha+\beta+2)_{2(k+m)+r-3}}\right\}^2
{}_4F_3
\left(
\begin{matrix}
-m,
-\frac{\beta-1}{2}-m,
\frac{\gamma}{2},
-\frac{\gamma}{2}
\cr
\frac{1}{2},
\frac{\alpha+r+1}{2},
-2m-\frac{\alpha+\beta+r-1}{2}
\end{matrix}; 1
\right)
\label{eq:cor-even}
\\&=(-2)^{3m}
\prod_{k=1}^{2m}\frac{(k-1)!(\alpha+1)_{k+r-1}}
{(\alpha+\beta+2)_{k+2m+r-2}}
\prod_{k=1}^{m}(\beta+1)_{2k-2}^2
\nonumber\\&\quad\times
W_{m}\left(-\frac{\gamma^2}4\,;\,0,\frac12,\frac{\alpha+r+1}2,-2m-\frac{\alpha+\beta+r-1}2\right).
\label{eq:cor-even2}
\end{align}
\item[(ii)]
If the size $n=2m+1$ of the matrix is odd,
then we have
\begin{align}
&\det\left((\gamma+j-i)\frac{(\alpha+1)_{i+j+r-2}}{(\alpha+\beta+2)_{i+j+r-2}}\right)_{1\leq i,j\leq 2m+1}\nonumber\\
&=\gamma
\prod_{k=1}^{m+1}\frac{(2k-1)!(\alpha+1)_{2k+r-2}(\beta+1)_{2k-2}}{(\alpha+\beta+2)_{2(k+m-1)+r}}
\prod_{k=1}^{m}\frac{(2k-1)!(\alpha+1)_{2k+r}(\beta+1)_{2k-2}}{(\alpha+\beta+2)_{2(k+m-1)+r}}
\nonumber\\&\quad\times
{}_4F_3
\left(
\begin{matrix}
-m,
-\frac{\beta-1}{2}-m, 
\frac{1+\gamma}{2},
\frac{1-\gamma}{2}
\cr
\frac{3}{2}, 
\frac{\alpha+r}{2}+1, 
-2m-\frac{\alpha+\beta+r}{2}
\end{matrix}; 1
\right)
\label{eq:cor-odd}
\\&=(-2)^{3m}\gamma
\prod_{k=1}^{2m+1}\frac{(k-1)!(\alpha+1)_{k+r-1}}
{(\alpha+\beta+2)_{k+2m+r-1}}
\prod_{k=1}^{m+1}(\beta+1)_{2k-2}
\prod_{k=1}^{m}(\beta+1)_{2k-2}
\nonumber\\&\quad\times
W_{m}\left(-\frac{\gamma^2}4\,;\,\frac12,1,\frac{\alpha+r+1}2,-2m-\frac{\alpha+\beta+r+1}2\right).
\label{eq:cor-odd2}
\end{align}
\end{enumerate}
\end{corollary}
%
%
%
%
%
%
%
%
This paper is composed as follows.
In Section~\ref{sec:main} we give a generalization of Theorem~\ref{th:corthm2} 
which is a determinant formula for arbitrary rows (see Theorem~\ref{lem:thm1}).
The most labors to prove our theorems exist in the proof of this formula.
In Section~\ref{sec:proof-theorems} we prove the main theorems in this section
from Theorem~\ref{lem:thm1}, which will be straightforward.
In Section~\ref{sec:quadratic-relation} we derive a quadratic relation among
the Askey-Wilson polynomials as a corollary of Theorem~\ref{th:corthm2}.
%
%
%
%
%
%
%
%
%
%
%
%
%
%

%% file: mehta02.tex
%
%
%
%
%
%
%
%
%
%
\section{Determinant formula for arbitrary rows}
\label{sec:main}
In our previous paper \cite{ITZ1},
we prove the following formula
in which the rows are arbitrary chosen.
Let $n$ be a positive integer,
and $k_{1}$, $\dots$, $k_{n}$ be arbitrary positive integers.
Then we have
\begin{align}
  &\det\left(
  \frac{(aq;q)_{k_{i}+j-2}}{(abq^2;q)_{k_{i}+j-2}}
  \right)_{1\leq i,j\leq n}
=a^{\frac{n(n-1)}2}q^{\frac{(n+1)n(n-1)}6}
\nonumber\\&\qquad\times
  \prod_{i=1}^{n}\frac{(aq;q)_{k_{i}-1}}{(abq^{2};q)_{k_{i}+n-2}}
  \prod_{1\leq i<j\leq n}(q^{k_{i}-1}-q^{k_{j}-1})
  \prod_{j=1}^{n}(bq;q)_{j-1}.
\label{eq:q-kratt}
\end{align}
This formula is a generalization of \eqref{eq:q-general-hankel}
and a special case is obtained in \cite[Theorem~3]{Kr2}. 
In this section we give this type formula, i.e.,
Theorem~\ref{lem:thm1},
which is crucial to prove Theorem~\ref{th:corthm2}.
\par\smallskip
%
%
%
%
%
First we fix some notation.
If $a$ and $b$ are integers,
we write $[a,b]=\left\{\,x\in\Z\,|\,a\leq x\leq b\,\right\}$.
We also write $[n]=[1,n]$ for short.
If $S$ is a finite set and $r$ a nonnegative integer,
let $\binom{S}{r}$ denote the set of all $r$-element subsets of $S$.
Let $A$ be an $m\times n$ matrix.
If $\boldsymbol{i}=(i_1,\dots,i_r)$ is an $r$-tuple of positive integers
and $\boldsymbol{j}=(j_1,\dots,j_s)$ is an $s$-tuple of positive integers,
then let $A^{\boldsymbol{i}}_{\boldsymbol{j}}=A^{i_1,\dots,i_r}_{j_1,\dots,j_s}$ denote the submatrix
formed by selecting the row $\boldsymbol{i}$ and the column $\boldsymbol{j}$ from $A$.
The main aim of this section is to prove the following key theorem.
%
%
%
%
\begin{theorem}
\label{lem:thm1}
Let $a$, $b$ and $c$ be parameters.
Let $n$ be a positive integer, and $\k=(k_{1},\dots,k_{n})$ be an $n$-tuple of positive integers.
Then we have
\begin{align}
&\det\left(
(q^{k_i-1}-cq^{j-1})
\frac{(aq;q)_{k_i+j-2}}{(abq^2;q)_{k_i+j-2}}
\right)_{1\leq i,j\leq n}\nonumber\\
&=
a^{\frac{n(n-3)}{2}}q^{\frac{n(n+1)(n-4)}{6}}
\prod_{i=1}^n\frac{(aq;q)_{k_i-1}(bq;q)_{i-2}}{(abq^2;q)_{k_i+n-2}}
\prod_{1\leq i<j\leq n}(q^{k_i-1}-q^{k_j-1})
\nonumber\\&\quad\times
\sum_{\nu=0}^n
(-1)^{n-\nu}(abcq^{2\nu+1};q^2)_{n-\nu}(acq;q^2)_\nu R_{n,\nu}(\k,a,b;q),
\label{eqextttqab}
\end{align}
where
\begin{align}
&R_{n,\nu}(\k,a,b;q)=
\sum_
{(\i,\j)}
q^{\sum_{l=1}^{n-\nu}i_l-n+\nu}
\nonumber\\&\qquad\qquad\times
\prod_{l=1}^{n-\nu}(1-aq^{k_{i_l}-i_l+l+\nu})
\prod_{l=1}^\nu(1-abq^{k_{j_l}+j_l-l+\nu-1}).
\label{eq:R_n}
\end{align}
Here the sum on the right-hand side runs over 
all pairs $(\i,\j)$
such that $[n]$ is a disjoint union of 
$\i=\{i_{1},\dots,i_{n-\nu}\}\in\binom{[n]}{n-\nu}$ and
$\j=\{j_{1},\dots,j_{\nu}\}\in\binom{[n]}{\nu}$
 (i.e., $\i\cup\j=[n]$ and $\i\cap\j=\emptyset$).
Hereafter, we also use the convention that $R_{n,\nu}(\k,a,b;q)$ is $0$ unless $0\leq\nu\leq n$.
\end{theorem}
For example,
if $n=3$ and $\nu=2$,
then the pairs $(\i,\j)$ runs over
\[
\left\{
\left(\{1\},\{2,3\}\right),(\{2\},\{1,3\}),(\{3\},\{1,2\})
\right\}.
\]
Hence we have
\begin{align*}
&R_{3,2}(\{k_1,k_2,k_3\},a,b;q)
=(1-aq^{k_1+2})(1-abq^{k_2+2})(1-abq^{k_3+2})
\\&\qquad\qquad
+q(1-aq^{k_2+1})(1-abq^{k_1+1})(1-abq^{k_3+2})
\\&\qquad\qquad
+q^2(1-aq^{k_3})(1-abq^{k_1+1})(1-abq^{k_2+1}).
\end{align*}
%
%
%
%
%
%
Let $n$ be a positive integer, and let $a$, $b$, $c$ and $q$ be parameters.
For an index set $\k=\{k_{1},\dots,k_{n}\}$ of positive integers,
let $M_{n}(\k,a,b,c;q)=\left(M_{n}(\k,a,b,c;q)_{i,j}\right)_{1\leq i,j\leq n}$ 
denote the matrix whose $(i,j)$ entry is given by
\begin{align}
\label{eq:M_n}
M_{n}(\k,a,b,c;q)_{i,j}=(q^{k_i-1}-cq^{j-1})\,(aq^{k_i};q)_{j-1}(abq^{k_i+j};q)_{n-j}.
\end{align}
Then we have
\begin{align}
\label{eq:AM_n}
\det\left(
(q^{k_i-1}-cq^{j-1})
\frac{(aq;q)_{k_i+j-2}}{(abq^2;q)_{k_i+j-2}}
\right)_{1\leq i,j\leq n}
=\prod_{i=1}^{n}\frac{(aq;q)_{k_i-1}}{(abq^2;q)_{k_i+n-2}}\cdot
\det M_{n}(\k,a,b,c;q).
\end{align}
Hence it is enough to evaluate $\det M_{n}(\k,a,b,c;q)$ to prove Theorem~\ref{lem:thm1}.
The main task of this evaluation is to show the following recurrence equation:
\begin{align}
&
\frac{\det M_{n}(\k,a,b,c;q)}
{a^{n-2}(bq;q)_{n-2}\prod_{i=1}^{n-1}(q^{k_i}-q^{k_n})}
\nonumber\\&
=
q^{-1}(1-acq)(1-abq^{k_n+n-1})\det M_{n-1}(\k',aq,b,cq;q)
\nonumber\\&\qquad\qquad
-q^{n(n-3)/2}(1-abcq^{2n-1})(1-aq^{k_n})\det M_{n-1}(\k',a,b,c;q),
\label{eq:rec_M_n}
\end{align}
where $\k'=\{k_1,\dots,k_{n-1}\}$ denote the subset of the first $(n-1)$ indices
of $\k=\{k_1,\dots,k_{n-1},k_{n}\}$.
This identity enable us to prove Theorem~\ref{lem:thm1} by induction.
%
%
%
%
First we cite the $q$-binomial formula \cite{GR,KLS,KS,Sta} (see also \cite[Lemma~2]{Ni}).
%
%
%
\begin{prop}
Let $n$ be a nonnegative integer.
Then we have
\begin{align}
&\sum_{k=0}^{n}
(-1)^{k}x^{k}q^{k(k-1)/2}\qbinom{n}{k}
=(x;q)_{n}.
\label{eq:Chu-Vander}
\end{align}
\end{prop}
%
%
%
%
%
%
\par\smallskip
Our method to evaluate $\det M_{n}(\k,a,b,c;q)$ is completely different from Mehta-Wang's proof or
Nishizawa's $q$-analogue.
It seems very tough to generalize their proof to our case.
So we establish the inductive identity \eqref{eq:rec_M_n} and appeal to the induction on the matrix size.
For the purpose
the following proposition plays an essential role in the matrix multiplication of $M_{n}(\k,a,b,c;q)$
which we use in the proof of Lemma~\ref{lem:First_Lemma}.
%
%
%
%
%
%
%
\begin{lemma}
\label{prop:residue}
Let $n$ be a positive integer.
Let $a$, $b$, $c$ and $q$ be complex numbers,
and $x_1$,$\hdots$, $x_n$ be variables.
Then we have
\begin{align}
&
-\sum_{\nu=1}^{n}
\frac{
(q^{-1}x_\nu-cq^{j-1})
(a x_\nu;q)_{j-1}(abq^{j}x_\nu;q)_{n-j}
}
{
x_\nu(1-ax_\nu)\prod_{{l=1}\atop{l\neq \nu}}^{n}(x_l-x_\nu)
}
\nonumber\\&\qquad
=
\frac{cq^{j-1}}{\prod_{l=1}^{n}x_l}
+\chi(j=1)\cdot\frac{(-1)^na^{n-1}q^{-1}(1-acq)(bq;q)_{n-1}}{\prod_{l=1}^{n}(1-ax_l)},
\label{eq:prop01}
\\&
-\sum_{\nu=1}^{n}
\frac{
(q^{-1}x_\nu-cq^{j-1})
(ax_\nu;q)_{j-1}(abq^{j}x_\nu;q)_{n-j}
}
{
x_\nu(1-abq^{n-1}x_\nu)\prod_{{l=1}\atop{l\neq \nu}}^{n}(x_l-x_\nu)
}
\nonumber\\&\qquad
=
\frac{cq^{j-1}}{\prod_{l=1}^{n}x_l}
-\chi(j=n)\cdot\frac{
a^{n-1}q^{n(n-3)/2}(1-abcq^{2n-1})
(bq;q)_{n-1}
}{\prod_{l=1}^{n}(1-abq^{n-1}x_l)}.
\label{eq:prop02}
\end{align}
\end{lemma}
%
%
%
%
\begin{demo}{Proof}
There are several methods to prove these identities.
Here we present a simple proof by the complex analysis
(see \cite{IKO}).
Let $F(z)$ be the meromorphic function of $z$ defined by 
\[
F(z)=\frac{
(q^{-1}z-cq^{j-1})
(az;q)_{j-1}(abq^{j}z;q)_{n-j}
}
{
z(1-az)\prod_{l=1}^{n}(x_l-z)
}.
\]
Then $z=x_\nu$ ($\nu=1,\dots,n$) is a pole, and its residue is
\begin{align*}
\underset{z=x_\nu}{\Res}\, F(z)=-\frac{
(q^{-1}x_\nu-cq^{j-1})(ax_\nu;q)_{j-1}(abq^{j}x_\nu;q)_{n-j}
}
{
x_\nu(1-ax_\nu)\prod_{{l=1}\atop{l\neq\nu}}^{n}(x_l-x_\nu)
}.
\end{align*}
Similarly $z=0$ is also a pole, and its residue is
\begin{align*}
\underset{z=0}{\Res}\,F(z)=-\frac{cq^{j-1}}
{\prod_{l=1}^{n}x_l}.
\end{align*}
At $z=a^{-1}$, $F(z)$ has the residue of
\begin{align*}
\underset{z=a^{-1}}{\Res}F(z)
&=-\frac{
(a^{-1}q^{-1}-cq^{j-1})
(1;q)_{j-1}(bq^{j};q)_{n-j}
}
{
\prod_{l=1}^{n}(x_l-a^{-1})
}\\
&=
-\chi(j=1)\cdot\frac{(a^{-1}q^{-1}-c)(bq;q)_{n-1}}{\prod_{l=1}^{n}(x_l-a^{-1})}.
\end{align*}
Finally $z=\infty$ is also a pole of $F(z)$, and
$
\underset{z=\infty}{\Res}F(z)=-\lim\limits_{z\rightarrow\infty} zF(z)=0.
$
Since the sum of the residues 
of a meromorphic function on a compact Riemann surface must be zero,
we obtain the desired identity \eqref{eq:prop01}.
The other identities can be proven similarly. 
The details are left to the reader. 
%
%
%
%
%
%
%
%
%
%
\end{demo}
Notice that an immediate consequence of these identities are the following Vandermonde type determinants.
These identities are obtained from Lemma~\ref{prop:residue}
by expanding the determinants along the last columns.
%
%
%
\begin{corollary}
Let $n$ be a positive integer, and $k$ be an integer such that $1\leq k\leq n$.
Let $a$, $b$, $c$ and $q$ be parameters,
and $\x=(x_1,\hdots,x_n)$ be an $n$-tuple of variables.
\begin{enumerate}
\item[(i)]
Let $V_{n,k}(\x,a,b,c;q)=\left(V_{n,k}(\x,a,b,c;q)_{i,j}\right)_{1\leq i,j\leq n}$ be the $n\times n$ matrix
defined by
\[
V_{n,k}(\x,a,b,c;q)_{i,j}
=\begin{cases}
x_{i}^{j-1}&\text{ if $1\leq j<n$,}\\
-\frac{
(x_i-cq^{k})
(a x_i;q)_{k-1}(abq^{k}x_i;q)_{n-k}
}
{
x_i(1-ax_i)
}
&\text{ if $j=n$,}
\end{cases}
\]
for $1\leq i\leq n$.
Then we have
\begin{align*}
&\frac{(-1)^{n-1}\det V_{n,k}(\x,a,b,c;q)}{\prod_{1\leq i<j\leq n}(x_j-x_i)}
\\&\quad
=
\frac{cq^{k}}{\prod_{l=1}^{n}x_l}
+\chi(k=1)\cdot\frac{(-1)^na^{n-1}(1-acq)(bq;q)_{n-1}}{\prod_{l=1}^{n}(1-ax_l)}.
\end{align*}
\item[(ii)]
Let $W_{n,k}(\x,a,b,c;q)=\left(W_{n,k}(\x,a,b,c;q)_{i,j}\right)_{1\leq i,j\leq n}$ be the $n\times n$ matrix
defined by
\[
W_{n,k}(\x,a,b,c;q)_{i,j}
=\begin{cases}
x_{i}^{j-1}&\text{ if $1\leq j<n$,}\\
-\frac{
(x_i-cq^{k})
(ax_i;q)_{k-1}(abq^{k}x_i;q)_{n-k}
}
{
x_{i}(1-abq^{n-1}x_{i})
}
&\text{ if $j=n$,}
\end{cases}
\]
for $1\leq i\leq n$.
Then we have
\begin{align*}
&\frac{(-1)^{n-1}\det W_{n,k}(\x,a,b,c;q)}{\prod_{1\leq i<j\leq n}(x_j-x_i)}
\\&\quad
=
\frac{cq^{k}}{\prod_{l=1}^{n}x_l}
-\chi(k=n)\cdot\frac{
a^{n-1}q^{(n-1)(n-2)/2}(1-abcq^{2n-1})
(bq;q)_{n-1}
}{\prod_{l=1}^{n}(1-abq^{n-1}x_l)}.
\end{align*}
\end{enumerate}
\end{corollary}
Here we will not use these Vandermonde type determinants,
but it seems that these identities are worth mentioning.
We introduce four triangular matrices $X_{n}(\k,a;q)$,
$Y_n(q)$, $L_{n}(\k,a,b;q)$ and $U_n(q)$ which plays an important role
to manipulate $M_{n}(\k,a,b,c;q)$ in \eqref{eq:AM_n}.
%
%
%
%
Let $X_{n}(\k,a;q)=\left(X(\k,a;q)_{i,j}\right)_{1\leq i,j\leq n}$ and
$Y_n(q)=\left(Y_n(q)_{i,j}\right)_{1\leq i,j\leq n}$
be the $n\times n$ lower triangular matrices whose $(i,j)$-entry is, respectively, given by
\begin{align}
&X(\k,a;q)_{i,j}
=-\frac{\chi(i\geq j)}{q^{k_{j}}(1-aq^{k_{j}})\prod_{{l=1}\atop{l\neq j}}^{i}(q^{k_{l}}-q^{k_{j}})},
\label{eq:X_n}
\\
&Y_n(q)_{i,j}=
(-1)^{i+j}q^{-\frac{(i-j)(2n+1-i-j)}2}\qbin{n-j}{i-j}.
\label{eq:Y_n}
\end{align}
%
%
%
Similarly,
let $L_{n}(\k,a,b;q)=\left(L_{n}(\k,a,b;q)_{i,j}\right)_{1\leq i,j\leq n}$ 
(resp. $U_n(q)=\left(U(q)_{i,j}\right)_{1\leq i,j\leq n}$ )
be the $n\times n$ lower (resp. upper) triangular matrix 
whose $(i,j)$-entry is, respectively, given by
\begin{align}
&L_{n}(\k,a,b;q)_{i,j}
=-\frac{\chi(i\geq j)}{q^{k_{j}}(1-abq^{k_{j}+n-1})\prod_{{l=1}\atop{l\neq j}}^{i}(q^{k_{l}}-q^{k_{j}})},
\label{eq:L_n}
\\
&U(q)_{i,j}=
(-1)^{i+j}q^{\frac{(j-i)(j-i+1)}2}\qbin{j-1}{j-i}.
\label{eq:U_n}
\end{align}
The fact that these are all triangular matrices is very important in the following proof.
But, note that these triangular matrices are irrelevant to the $LU$-decomposition of $M_{n}(\k,a,b,c;q)$.
The $LU$-decomposition of $M_{n}(\k,a,b,c;q)$ seems another difficult problem.
We define the $n\times n$ matrices $P_{n}(\k,a,b,c;q)$ and $Q_{n}(\k,a,b,c;q)$ by
\begin{align*}
&
P_{n}(\k,a,b,c;q)=X_{n}(\k,a;q)M_{n}(\k,a,b,c;q)Y_n(q),
\\&
Q_{n}(\k,a,b,c;q)=L_{n}(\k,a,b;q)M_{n}(\k,a,b,c;q)U_n(q).
\end{align*}
Since $X_{n}(\k,a;q)$, $L_{n}(\k,a,b;q)$ are triangular
and $Y_{n}(q)$, $U_{n}(q)$ are unitriangular,
we easily obtain
\begin{align}
&
\det P_{n}(\k,a,b,c;q)=\frac{(-1)^n\det M_{n}(\k,a,b,c;q)}{q^{\sum_{i=1}^{n}k_i}\prod_{i=1}^{n}(1-aq^{k_i})\prod_{1\leq i<j\leq n}(q^{k_i}-q^{k_j})},
\label{eq:P_{n}M_{n}}
\\&
\det Q_{n}(\k,a,b,c;q)=\frac{(-1)^n\det M_{n}(\k,a,b,c;q)}{q^{\sum_{i=1}^{n}k_i}\prod_{i=1}^{n}(1-abq^{k_i+n-1})\prod_{1\leq i<j\leq n}(q^{k_i}-q^{k_j})}.
\label{eq:Q_{n}M_{n}}
\end{align}
The key to prove \eqref{eq:rec_M_n} will be Lemma~\ref{lem:PM-QM}. 
To prove Lemma~\ref{lem:PM-QM}, we need the following lemma
which gives the bottom rows of $P_{n}(\k,a,b,c;q)$ and $Q_{n}(\k,a,b,c;q)$.
%
%
%
%
\begin{lemma}
\label{lem:First_Lemma}
Let $n\geq2$ be an integer,
and let $a$, $b$, $c$ and $q$ be parameters.
Let $X_{n}(\k,a;q)$, $Y_n(q)$, $L_{n}(\k,a,b;q)$, $U_n(q)$, $P_{n}(\k,a,b,c;q)$ and $Q_{n}(\k,a,b,c;q)$ be as above,
and $M_{n}(\k,a,b,c;q)$ as in \eqref{eq:M_n}.
Then the bottom rows of $P_{n}(\k,a,b,c;q)$ and $Q_{n}(\k,a,b,c;q)$ are given by
\begin{align}
&
P_{n}(\k,a,b,c;q)_{n,j}
=\begin{cases}
\frac{(-1)^{n}a^{n-1}q^{-1}(1-acq)(bq;q)_{n-1}}{\prod_{l=1}^{n}(1-aq^{k_l})}
&\text{ if $j=1$,}\\
0&\text{ if $1<j<n$,}\\
cq^{n-1-\sum_{l=1}^{n}k_l}&\text{ if $j=n$,}
\end{cases}
\label{eq:lem01}
\\&
Q_{n}(\k,a,b,c;q)_{n,j}
=\begin{cases}
cq^{-\sum_{l=1}^{n}k_l}&\text{ if $j=1$,}\\
0&\text{ if $1< j<n$,}\\
-\frac{
a^{n-1}q^{n(n-3)/2}(1-abcq^{2n-1})(bq;q)_{n-1}
}{\prod_{l=1}^{n}(1-abq^{k_l+n-1})}
&\text{ if $j=n$.}
\end{cases}
\label{eq:lem02}
\end{align}
\end{lemma}
%
%
%
%
%
\begin{demo}{Proof}
By substituting $x_{l}=q^{k_l}$ for $l=1,\dots,n$ into \eqref{eq:prop01},
we see the $(n,j)$-entry of $X_{n}(\k,a;q)M_{n}(\k,a,b,c;q)$ equals
\[
cq^{j-1-\sum_{l=1}^{n}k_l}
+\chi(j=1)\cdot\frac{(-1)^na^{n-1}q^{-1}(1-acq)(bq;q)_{n-1}}{\prod_{l=1}^{n}(1-aq^{k_l})}.
\]
Since $Y_{n}(q)$ is lower unitriangular,
the $(n,1)$-entry of $X_{n}(\k,a;q)M_{n}(\k,a,b,c;q)$ affects only the first entry of the bottom row
of $X_{n}(\k,a;q)M_{n}(\k,a,b,c;q)Y_{n}(q)$.
Hence, if $j\neq1$, then we have
\[
P_{n}(\k,a,b,c;q)_{n,j}=\sum_{\nu=j}^{n}
(-1)^{\nu+j}cq^{\nu-1-\frac{(\nu-j)(2n+1-\nu-j)}2-\sum_{l=1}^{n}k_l}
\qbin{n-j}{\nu-j}.
\]
%
%
By replacing $\nu$ by $n-\nu$, we obtain
\[
P_{n}(\k,a,b,c;q)_{n,j}=cq^{j-1-\sum_{l=1}^{n}k_l-\frac{(n-j)(n-j-1)}2}
\sum_{\nu=0}^{n-j}
(-1)^{n-j-\nu}q^{\frac{\nu(\nu-1)}2}
\qbin{n-j}{\nu}.
\]
By \eqref{eq:Chu-Vander}, this equals
$cq^{n-1-\sum_{l=1}^{n}k_l}$ if $j=n$,
and $0$ otherwise.
The case where $j=1$ is easily obtained from $Y_{n}(q)_{1,1}=1$.
This proves \eqref{eq:lem01}.
The other identities can be proven similarly. 
The details are left to the reader. 
%
%
%
%
%
%
\end{demo}
%
%
%
%
%
\begin{prop}
\label{eq:minorsYU}
Let $n$ be a positive integer,
and $q$ a parameter.
Let $Y_{n}(q)$ and $U_{n}(q)$ be as in \eqref{eq:Y_n} and \eqref{eq:U_n}.
Then we have
\begin{align}
&
Y_{n}(q)^{-1}=\left(q^{(j-i)(n+1-i)}\qbinom{n-j}{i-j}\right)_{1\leq i,j\leq n},
\label{eq:inverse-Y}
\\&
U_{n}(q)^{-1}=\left(q^{j-i}\qbinom{j-1}{i-1}\right)_{1\leq i,j\leq n}.
\label{eq:inverse-U}
\end{align}
Especially,
we have
\begin{align}
&\det Y_{n}(q)^{[1,n]\setminus\{i\}}_{[1,n-1]}
=(-q)^{i-n},
\label{eq:sub_Y_n}\\
&\det U_{n}(q)^{[1,n]\setminus\{i\}}_{[2,n]}
=(-q)^{i-1}.
\label{eq:sub_U_n}
\end{align}
\end{prop}
%
%
%
%
%
%
\begin{demo}{Proof}
%
To prove that \eqref{eq:inverse-Y} gives the inverse of $Y_{n}(q)$,
one need to show that
\begin{align*}
&\sum_{k=j}^{i}q^{(k-i)(n+1-i)}\qbinom{n-k}{i-k}(-1)^{k+j}q^{-\frac{(k-j)(2n+1-k-j)}2}
\qbinom{n-j}{k-j}
\\&\quad=
(-1)^{i+j}q^{-\frac{(i-j)(2n+1-i-j)}2}\qbinom{n-j}{n-i}
\sum_{k=0}^{i-j}(-1)^{k}q^{\frac{k(k-1)}2}\qbinom{i-j}{k}
\end{align*}
equals $\chi(i=j)$,
where the second sum is obtained from the first sum by replacing $i-k$ by $k$.
This can be shown by \eqref{eq:Chu-Vander}.
It is also an easy exercise to show that \eqref{eq:inverse-U}
gives the inverse matrix of $U_{n}(q)$ using \eqref{eq:Chu-Vander}.
Finally,
\eqref{eq:sub_Y_n} (resp. \eqref{eq:sub_U_n}) is obtained from
\eqref{eq:inverse-Y} (resp. \eqref{eq:inverse-U})
using the relation $A^{-1}=\frac1{\det A}\adj(A)$,
where the adjugate matrix $\adj(A)$ is the transpose of the matrix of cofactors.
\end{demo}
Let $n\leq N$ be positive integers,
and let $A$ be an $n\times N$ matrix and $B$ an $N\times n$ matrix.
Then the following formula is known as the Cauchy--Binet formula:
\begin{equation}
\det AB
=\sum_{\boldsymbol{i}\in\binom{[N]}{n}} A^{[n]}_{\boldsymbol{i}}B^{\boldsymbol{i}}_{[n]}.
\label{eq:Cauchy-Binet}
\end{equation}
%
%
%
\begin{lemma}
\label{lem:PM-QM}
Let $n$ be a positive integer,
and let $a$, $b$, $c$ and $q$ be parameters.
Let $P_{n}(\k,a,b,c;q)$ and $Q_{n}(\k,a,b,c;q)$ be as defined in Lemma~\ref{lem:First_Lemma}.
When $\k=\{k_1,\dots,k_{n-1},k_{n}\}$ is a row index set,
let $\k'=\{k_1,\dots,k_{n-1}\}$ denote the subset of the first $(n-1)$ indices of $\k$.
Then we have
\begin{align}
&
\det P_{n}(\k,a,b,c;q)^{[1,n-1]}_{[2,n]}
=\frac{(-1)^{n-1}\det M_{n-1}(\k',aq,b,cq;q)}
{q^{\sum_{i=1}^{n-1}k_i}\prod_{1\leq i<j<n}(q^{k_i}-q^{k_j})},
\label{eq:P_n}
\\&
\det Q_{n}(\k,a,b,c;q)^{[1,n-1]}_{[1,n-1]}
=\frac{(-1)^{n-1}\det M_{n-1}(\k',a,b,c;q)}
{q^{\sum_{i=1}^{n-1}k_i}\prod_{1\leq i<j<n}(q^{k_i}-q^{k_j})},
\label{eq:Q_n}
\\&
\frac{\det P_{n}(\k,a,b,c;q)^{[1,n-1]}_{[1,n-1]}}
{\prod_{\nu=1}^{n-1}(1-abq^{k_\nu+n-1})}
=(-q)^{-n+1}\frac{\det Q_{n}(\k,a,b,c;q)^{[1,n-1]}_{[2,n]}}
{\prod_{\nu=1}^{n-1}(1-aq^{k_\nu})}.
\label{eq:P_nQ_n}
\end{align}
\end{lemma}
%
%
%
\begin{demo}{Proof}
Since $X_{n}(\k,a;q)$ and $Y_n(q)$ are lower triangular,
we have 
\[
P_{n}(\k,a,b,c;q)^{[1,n-1]}_{[2,n]}=X_{n}(\k,a;q)^{[1,n-1]}_{[1,n-1]}\,
M_{n}(\k,a,b,c;q)^{[1,n-1]}_{[2,n]}\,Y_n(q)^{[2,n]}_{[2,n]}.
\]
Hence,
\eqref{eq:P_n} easily follows from
\begin{align*}
&
\det X_{n}(\k,a;q)^{[1,n-1]}_{[1,n-1]}
=\frac{(-1)^{n-1}}
{q^{\sum_{\nu=1}^{n-1}k_\nu}\prod_{\nu=1}^{n-1}(1-aq^{k_\nu})\prod_{1\leq i<j<n}(q^{k_i}-q^{k_j})},
\\&
\det M_{n}(\k,a,b,c;q)^{[1,n-1]}_{[2,n]}
=\prod_{\nu=1}^{n-1}(1-aq^{k_\nu})\cdot\det M_{n-1}(\k',aq,b,cq;q),
\\&
\det Y_n(q)^{[2,n]}_{[2,n]}=1.
\end{align*}
Here, the first and the third identities follow from the fact that
$X_{n}(\k,a;q)$ and $Y_n(q)$ are lower triangular,
and the second identity follows from the fact that 
$(i,j)$ entry of $M_{n}(\k,a,b,c;q)^{[1,n-1]}_{[2,n]}$
equals $(q^{k_i-1}-cq^{j})(aq^{k_i};q)_{j}(abq^{k_i+j+1};q)_{n-j-1}$ for $1\leq i,j<n$
(see \eqref{eq:M_n}).
Exactly the same argument proves \eqref{eq:Q_n} from the fact
\[
Q_{n}(\k,a,b,c;q)^{[1,n-1]}_{[1,n-1]}=L_{n}(\k,a,b;q)^{[1,n-1]}_{[1,n-1]}
\,M_{n}(\k,a,b,c;q)^{[1,n-1]}_{[1,n-1]}\,U_n(q)^{[1,n-1]}_{[1,n-1]}.
\]
The details of this argument are left to the reader.
%
Finally,
we prove \eqref{eq:P_nQ_n}.
If there is no fear of confusion,
we may write $M_{n}(\k,a,b,c;q)$ (resp. 
$P_{n}(\k,a,b,c;q)$, $Q_{n}(\k,a,b,c;q)$, $X_{n}(\k,a;q)$, $Y_n(q)$,
$L_{n}(\k,a,b;q)$ and $U_n(q)$) 
as $M_{n}(\k)$  (resp. $P_{n}(\k)$, $Q_{n}(\k)$, $X_{n}(\k)$, $Y_n$,
$L_{n}(\k)$ and $U_n$) in short hereafter.
We use the identities
\begin{align*}
&P_{n}(\k)^{[1,n-1]}_{[1,n-1]}=X_{n}(\k)^{[1,n-1]}_{[1,n-1]}\,
M_{n}(\k)^{[1,n-1]}_{[1,n]}\,{Y_n}^{[1,n]}_{[1,n-1]},
\\&
Q_{n}(\k)^{[1,n-1]}_{[2,n]}=L_{n}(\k)^{[1,n-1]}_{[1,n-1]}
\,M_{n}(\k)^{[1,n-1]}_{[1,n]}\,{U_n}^{[1,n]}_{[2,n]},
\end{align*}
which come from the triangularities of the matrices as before.
Hence,
by taking the determinant of the both sides,
we obtain
\begin{align*}
&
\det P_{n}(\k)^{[1,n-1]}_{[1,n-1]}=\det X_{n}(\k)^{[1,n-1]}_{[1,n-1]}\cdot\det M_{n}(\k)^{[1,n-1]}_{[1,n]}\,{Y_n}^{[1,n]}_{[1,n-1]},
\\&
\det Q_{n}(\k)^{[1,n-1]}_{[2,n]}=\det L_{n}(\k)^{[1,n-1]}_{[1,n-1]}\cdot\det M_{n}(\k)^{[1,n-1]}_{[1,n]}\,{U_n}^{[1,n]}_{[2,n]}.
\end{align*}
Since the $X_{n}(\k,a;q)^{[1,n-1]}_{[1,n-1]}$ and $L_{n}(\k,a,b;q)^{[1,n-1]}_{[1,n-1]}$ are lower triangular,
it is easy to compute the first determinant of the right-hand side of each equality.
Applying the Cauchy--Binet formula \eqref{eq:Cauchy-Binet} to the second determinant of the right-hand side of each equality,
we obtain
\begin{align*}
&
\det P_{n}(\k)^{[1,n-1]}_{[1,n-1]}
=\frac{(-1)^{n-1}
\sum_{i=1}^{n}
\det M_{n}(\k)^{[1,n-1]}_{[1,n]\setminus\{i\}} 
\det {Y_{n}}^{[1,n]\setminus\{i\}}_{[1,n-1]}
}{
q^{\sum_{\nu=1}^{n-1}k_\nu}\prod_{\nu=1}^{n-1}(1-aq^{k_\nu})\prod_{1\leq i<j<n}(q^{k_i}-q^{k_j})
},
\\&
\det Q_{n}(\k)^{[1,n-1]}_{[2,n]}
=\frac{(-1)^{n-1}
\sum_{i=1}^{n}
\det M_{n}(\k)^{[1,n-1]}_{[1,n]\setminus\{i\}}
\det {U_{n}}^{[1,n]\setminus\{i\}}_{[2,n]}
}{
q^{\sum_{\nu=1}^{n-1}k_\nu}\prod_{\nu=1}^{n-1}(1-abq^{k_\nu+n-1})\prod_{1\leq i<j<n}(q^{k_i}-q^{k_j})
}.
\end{align*}
Hence \eqref{eq:P_nQ_n} immediately follows from Proposition~\ref{eq:minorsYU}.
This completes the proof.
\end{demo}
Now we are in position to prove \eqref{eq:rec_M_n}.
%
%
%
\begin{demo}{Proof of \eqref{eq:rec_M_n}}
Expanding $P_{n}(\k)$ and $Q_{n}(\k)$ along the bottom row,
we obtain
\begin{align}
&
\det P_{n}(\k)
=-\frac{a^{n-1}q^{-1}(1-acq)(bq;q)_{n-1}}{\prod_{l=1}^{n}(1-aq^{k_l})}
\det {P_{n}(\k)}^{[1,n-1]}_{[2,n]}
\nonumber\\&\qquad\qquad\qquad\qquad
+cq^{n-1-\sum_{l=1}^{n}k_l}
\det {P_{n}(\k)}^{[1,n-1]}_{[1,n-1]},
\label{eq:{P_{n}}^{[1,n-1]}_{[1,n-1]}}
\\&
\det Q_{n}(\k)
=(-1)^{n+1}cq^{-\sum_{l=1}^{n}k_l}
\det {Q_{n}(\k)}^{[1,n-1]}_{[2,n]}
\nonumber\\&\qquad\qquad
-\frac{
a^{n-1}q^{n(n-3)/2}(1-abcq^{2n-1})(bq;q)_{n-1}
}{\prod_{l=1}^{n}(1-abq^{k_l+n-1})}
\det {Q_{n}(\k)}^{[1,n-1]}_{[1,n-1]},
\label{eq:{Q_{n}}^{[1,n-1]}_{[1,n-1]}}
\end{align}
from \eqref{eq:lem01} and \eqref{eq:lem02}.
Hence, substituting \eqref{eq:P_{n}M_{n}} and \eqref{eq:P_n} into \eqref{eq:{P_{n}}^{[1,n-1]}_{[1,n-1]}},
\eqref{eq:Q_{n}M_{n}} and \eqref{eq:Q_n} into \eqref{eq:{Q_{n}}^{[1,n-1]}_{[1,n-1]}},
then the resulting equalities into \eqref{eq:P_nQ_n},
we obtain the desired identity.
\end{demo}
Now we are in position to prove Theorem~\ref{lem:thm1}.
In fact the proof is straightforward by induction.
\begin{demo}{Proof of Theorem~\ref{lem:thm1}}
First, we note that, for any integers $n$ and $\nu$,
it holds
\begin{align}
R_{n,\nu}(\k,a,b;q)
&=(1-abq^{k_{n}+n-1})R_{n-1,\nu-1}(\k',aq,b;q)
\nonumber\\&\qquad\qquad
+q^{n-1}(1-aq^{k_n})R_{n-1,\nu}(\k',a,b;q),
\label{eq:recRn_nu}
\end{align}
where $\k=\{k_1,\dots,k_{n-1},k_n\}$ and $\k'=\{k_1,\dots,k_{n-1}\}$
are as before.
\eqref{eq:recRn_nu} follows from the definition \eqref{eq:R_n} of $R_{n,\nu}(\k,a,b;q)$
by considering two exclusive cases, $j_\nu=n$ or $i_{n-\nu}=n$.
Now we prove the identity
\begin{align}
&\det M_{n}(\k,a,b,c;q)
=
(-1)^na^{\frac{n(n-3)}{2}}q^{\frac{n(n+1)(n-4)}{6}}
\prod_{i=1}^n (bq;q)_{i-2}
\nonumber\\&\quad\times
\prod_{1\leq i<j\leq n}(q^{k_i-1}-q^{k_j-1})
\sum_{\nu=0}^n
(-1)^\nu(abcq^{2\nu+1};q^2)_{n-\nu}(acq;q^2)_\nu\, R_{n,\nu}(\k,a,b;q),
\label{eq:induc_M_{n}}
\end{align}
by induction on $n$.
If $n=1$, then the left-hand side of \eqref{eq:induc_M_{n}} is trivially $q^{k_1-1}-c$
from \eqref{eq:M_n}.
It is straightforward computation to check the right-hand side equals $q^{k_1-1}-c$.
Assume $n>1$ and \eqref{eq:induc_M_{n}} holds up to $(n-1)$.
Using \eqref{eq:rec_M_n} and the induction hypothesis,
we obtain
\begin{align*}
&
\frac{\det M_{n}(\k,a,b,c;q)}
{(-1)^na^{\frac{n(n-3)}2}q^{\frac{n(n^2-6n-1)}6}\prod_{i=1}^{n}(bq;q)_{i-2}\prod_{1\leq i<j\leq n}(q^{k_i}-q^{k_j})}
\\&
=
(1-abq^{k_n+n-1})\sum_{\nu=0}^{n-1}(-1)^{\nu+1}(abcq^{2\nu+3};q^2)_{n-\nu-1}(acq;q^2)_{\nu+1}R_{n-1,\nu}(\k',aq,b;q)
\\&
+q^{n-1}(1-aq^{k_n})
\sum_{\nu=0}^{n-1}(-1)^{\nu}(abcq^{2\nu+1};q^2)_{n-\nu}(acq;q^2)_{\nu}R_{n-1,\nu}(\k',a,b;q).
\end{align*}
Replacing $\nu+1$ by $\nu$ in the first sum
and applying \eqref{eq:recRn_nu},
we establish \eqref{eq:induc_M_{n}} for $n$.
Hence \eqref{eq:induc_M_{n}} holds for an
arbitrary positive integer $n$.
Finally, \eqref{eq:AM_n} and \eqref{eq:induc_M_{n}} 
immediately implies \eqref{eqextttqab}.
This completes the proof of Theorem~\ref{lem:thm1}.
\end{demo}
To recover \eqref{eq:q-kratt} from Theorem~\ref{lem:thm1},
one can prove
\begin{equation}
\sum_{\nu=0}^n(-1)^{n-\nu}R_{n,\nu}({\boldsymbol{k}},a,b;q)
=a^nq^{\frac{n(n-1)}2+\sum_{l=1}^nk_l}(b;q)_n
\label{eq:sumRn}
\end{equation}
by induction, appealing to \eqref{eq:recRn_nu}.
If he substitutes \eqref{eq:sumRn} into \eqref{eqextttqab},
then he obtains \eqref{eq:q-kratt} immediately.
%
%
%
%
%
%
%
%
%
%
%
%
%
%
%
%
%
%
%
%
%
%
%
%
%
%
%
%
%
%
%
%
%
%
%

%% file: mehta03.tex

\section{Proof of the main theorems}
\label{sec:proof-theorems}

The aim of this section is 
to derive Theorem~\ref{th:corthm2} from Theorem~\ref{lem:thm1},
and then prove Corollary~\ref{th:main} from Theorem~\ref{th:corthm2}.
Once we prove Theorem~\ref{lem:thm1},
then it is easy and straightforward to prove the main theorems mainly by induction.
First,
to prove Theorem~\ref{th:corthm2},
we set 
$\k=[n]=\{1,2,\dots,n\}$
in \eqref{eqextttqab}.
Hence, it is essential to prove the following lemma
before proceeding to the proof of \eqref{eqtnabcn3aw}.
%
%
%
%
\begin{lemma}
\label{lem:R-[n]}
If we put $\k=[n]$ in \eqref{eqextttqab},
then we obtain
\begin{align}
R_{n,\nu}([n],a,b;q)=
q^{\frac{(n-\nu)(n-\nu-1)}{2}}\left[{{n}\atop{\nu}}\right]_q(aq^{\nu+1};q)_{n-\nu}(abq^n;q)_\nu.
\label{eqr123abc}
\end{align}
\end{lemma}
%
%
%
%
\begin{demo}{Proof}
We proceed by induction on $n$.
If $n=0$, then the both-sides are equal to $1$ if $\nu=0$, and $0$ otherwise.
Let $n>0$,
and assume \eqref{eqr123abc} holds when $n-1$.
By \eqref{eq:recRn_nu},
the left-hand side satisfies
\begin{align*}
R_{n,\nu}([n],a,b;q)
&=(1-abq^{2n-1})R_{n-1,\nu-1}([n-1],aq,b;q)
\\&\qquad\qquad
+q^{n-1}(1-aq^{n})R_{n-1,\nu}([n-1],a,b;q).
\end{align*}
Substituting the induction hypothesis into this identity,
we see the right-hand side equals
\begin{align*}
&
q^{\frac{(n-\nu)(n-\nu-1)}2}(aq^{\nu+1};q)_{n-\nu}(abq^{n};q)_{\nu-1}
\\&\qquad\times
\biggl\{\qbinom{n-1}{\nu-1}(1-abq^{2n-1})+q^{\nu}\qbinom{n-1}{\nu}(1-abq^{n-1})\biggr\},
\end{align*}
which becomes the right-hand side of \eqref{eqr123abc}.
Hence this proves that \eqref{eqr123abc} holds for an arbitrary nonnegative number $n$.
\end{demo}
Now we are in position to prove Theorem~\ref{th:corthm2}
%
%
%
%
\begin{demo}{Proof of Theorem~\ref{th:corthm2}}
If we substitute $\k=[n]$ into \eqref{eqextttqab} using \eqref{eqr123abc},
then we see that
$
\det\left(
(q^{i-1}-cq^{j-1})\frac{(aq;q)_{i+j-2}}{(abq^2;q)_{i+j-2}}
\right)_{1\leq i,j\leq n}
$
equals
\begin{align*}
&
a^{\frac{n(n-3)}2}q^{\frac{n(n+1)(n-4)}6}
\prod_{k=1}^{n}\frac{(aq;q)_{k-1}(bq;q)_{k-2}}{(abq^2;q)_{k+n-2}}
\prod_{1\leq i<j\leq n}(q^{i-1}-q^{j-1})
\\&\times
\sum_{\nu=0}^{n}(-1)^{n-\nu}q^{\frac{(n-\nu)(n-\nu-1)}2}
\qbinom{n}{\nu}(abcq^{2\nu+1};q^2)_{n-\nu}
(acq;q^2)_{\nu}(aq^{\nu+1};q)_{n-\nu}(abq^{n};q)_{\nu}.
\end{align*}
By rewriting the sum with the basic hypergeometric series
by using $\prod_{1\leq i<j\leq n}(q^{i-1}-q^{j-1})=q^{\frac{n(n-1)(n-2)}6}\prod_{k=0}^{n-1}(q;q)_k$,
$\qbinom{n}{\nu}=(-1)^{\nu}q^{n\nu-\frac{\nu(\nu-1)}2}\frac{(q^{-n};q)_{\nu}}{(q;q)_{\nu}}$,
$(abcq^{2\nu+1};q^2)_{n-\nu}=\frac{(abcq;q^2)_{n}}{(abcq;q^2)_{\nu}}$
and
$(aq^{\nu+1};q)_{n-\nu}=\frac{(aq;q)_{n}}{(aq;q)_{\nu}}$,
this becomes
\begin{align*}
&
(-1)^{n}a^{\frac{n(n-3)}2}q^{\frac{n(n+1)(2n-5)}6}
(abcq;q^2)_{n}
\prod_{k=1}^{n}\frac{(q;q)_{k-1}(aq;q)_{k}(bq;q)_{k-2}}{(abq^2;q)_{k+n-2}}
\\&\qquad\qquad\times
{}_{4}\phi_{3}\left(
{
{a^{\frac12}c^{\frac12}q^{\frac12},-a^{\frac12}c^{\frac12}q^{\frac12},q^{-n},abq^{n}}
\atop
{a^{\frac12}b^{\frac12}c^{\frac12}q^{\frac12},-a^{\frac12}b^{\frac12}c^{\frac12}q^{\frac12},aq}
}\,;\,q,q\right).
\end{align*}
If we replace $a$ by $aq^{r}$, and multiply the both sides
with $\left(\frac{(aq;q)_{r}}{(abq^2;q)_{r}}\right)^n$,
then we obtain \eqref{eqtnabcn3}.
There are several ways to write \eqref{eqtnabcn3} with the Askey-Wilson polynomials,
and \eqref{eqtnabcn3aw} follows from \eqref{eqtnabcn3} using the definition \eqref{eq:Askey-Wilson}.
\end{demo}
Before we proceed to a proof of Corollary~\ref{th:main},
we need the following contiguous relations for ${}_{4}\phi_{3}$.
%
%
\begin{prop}
\label{prop:phi-identities}
Let $z$, $a$, $b$, $c$, $d$, $e$, $f$, $g$ and $q$ be arbitrary parameters.
Then we have
\begin{align}
%
%
&
{}_{4}\phi_{3}\left( {
{a,bq,c,d}
\atop
{e,f,g}
};q,z \right) 
-{}_{4}\phi_{3}\left( {
{aq,b,c,d}
\atop
{e,f,g}
};q,z\right) 
\nonumber\\&\qquad=
\frac{z\left( b-a \right)  ( 1-c )  ( 1-d ) }
{ ( 1-e ) ( 1-f ) ( 1-g ) }
{}_{4}\phi_{3} \left( {
{aq,bq,cq,dq}
\atop
{eq,fq,gq}
} ; q,z \right), 
\label{eq:phi01}
\\
%
%
&
(1-f)(a-e)\,{}_{4}\phi_{3}\left( {
{a,b,c,d}
\atop
{eq,f,g}
};q,z \right) 
-(1-e)(a-f)\,{}_{4}\phi_{3}\left( {
{a,b,c,d}
\atop
{e,fq,g}
};q,z\right) 
\nonumber\\&\qquad=
( 1-a ) ( f-e )  \,
{}_{4}\phi_{3} \left( {
{aq,b,c,d}
\atop
{eq,fq,g}
} ; q,z \right),
\label{eq:phi02}
\end{align}
and
\begin{align}
%
%
&
(1-e)(1-f)(1-g)
\,
{}_{4}\phi_{3} \left( {
{a,b,c,d}
\atop
{e,f,g}
} ; q,q \right)
\nonumber\\&=
c(1-e)\left(1-\frac{f}{c}\right)\left(1-\frac{g}{c}\right)
\,
{}_{4}\phi_{3} \left( {
{aq,bq,c,d}
\atop
{e,fq,gq}
} ; q,q \right)
\nonumber\\&+
d(1-c)\left(1-\frac{e}{d}\right)\left(1-\frac{fg}{cd}\right)
\,
{}_{4}\phi_{3} \left( {
{aq,bq,cq,d}
\atop
{eq,fq,gq}
} ; q,q \right),
\label{eq:phi04}
\end{align}
where, in the last identity, we assume $abcdq=efg$ and $a=q^{-n}$ for some nonnegative integer $n$.
\end{prop}
%
%
%
%
%
\begin{demo}{Proof}
First, \eqref{eq:phi01} and \eqref{eq:phi02} are readily proven by direct computation.
The last identity \eqref{eq:phi04} is written as
\begin{align}
&
(1-e)(1-f)(1-g)
\,
{}_{4}\phi_{3} \left( {
{q^{-n},\frac{efgq^{n-1}}{cd},c,d}
\atop
{e,f,g}
} ; q,q \right)
\nonumber\\&=
c(1-e)\left(1-\frac{f}{c}\right)\left(1-\frac{g}{c}\right)
\,
{}_{4}\phi_{3} \left( {
{q^{-n+1},\frac{efgq^{n}}{cd},c,d}
\atop
{e,fq,gq}
} ; q,q \right)
\nonumber\\&+
d(1-c)\left(1-\frac{e}{d}\right)\left(1-\frac{fg}{cd}\right)
\,
{}_{4}\phi_{3} \left( {
{q^{-n+1},\frac{efgq^{n}}{cd},cq,d}
\atop
{eq,fq,gq}
} ; q,q \right),
\label{eq:phi04b}
\end{align}
and can be proven by induction on $n$ by using only \eqref{eq:phi01}.
The details are left to the reader.
\end{demo}
%
%
%
%
%
%
%
%
%
\begin{remark}
\label{rem:phi-identities}
The contiguous relations \eqref{eq:phi01} (resp. \eqref{eq:phi02}) correspond to 
(3.2) (resp. (3.10)) in \cite{Kr3},
meanwhile \eqref{eq:phi02} can be written as a contiguous relation for ${}_{8}W_{7}$.
In fact, if one uses Watson's transformation formula \cite[(2.5.1)]{GR}
\begin{align}
&
{}_{8}W_{7}\left(a;b,c,d,e,q^{-n};q,\frac{a^2q^{n+2}}{bcde}\right)
=\frac{\left(aq,\frac{aq}{de};q\right)_{n}}{\left(\frac{aq}{d},\frac{aq}{e};q\right)_{n}}
\,
{}_{4}\phi_{3} \left( {
{q^{-n},d,e,\frac{aq}{bc}}
\atop
{\frac{aq}{b},\frac{aq}{c},\frac{deq^{-n}}{a}}
} ; q,q \right)
\end{align}
for a terminating very-well-poised ${}_{8}\phi_{7}$ series,
where
\begin{eqnarray}
&{}_{r+1}W_{r}(a_{1};a_{4},\dots,a_{r+1};q,z)
={}_{r+1}\phi_{r}\left(
{{a_{1},qa_{1}^{\frac12},-qa_{1}^{\frac12},a_{4},\dots,a_{r+1}}
\atop{a_{1}^{\frac12},-a_{1}^{\frac12},\frac{qa_{1}}{a_{4}},\dots,\frac{qa_{1}}{a_{r+1}}}}
\,;\,q,z
\right),
\end{eqnarray}
then \eqref{eq:phi04} is equivalent to
\begin{align}
&
(c-a)(d-aq)(e-aq)(b-aq^{n})\,
{}_{8}W_{7}\left(a;b,cq,d,e,q^{-n};q,\frac{a^2q^{n+1}}{bcde}\right)
\nonumber\\&=
a(1-b)(1-aq)(de-aq)(1-cq^{n})\,
{}_{8}W_{7}\left(aq;bq,cq,d,e,q^{-n+1};q,\frac{a^2q^{n+1}}{bcde}\right)
\nonumber\\&+
(bc-a)(d-aq)(e-aq)(1-aq^{n})\,
{}_{8}W_{7}\left(a;b,c,d,e,q^{-n+1};q,\frac{a^2q^{n+1}}{bcde}\right).
\end{align}
\end{remark}
The following proposition is crucial to derive Corollary~\ref{th:main}
from Theorem~\ref{th:corthm2}.
%
%
%
%
%
\begin{prop}
\label{lem:even-odd}
Let $n$ be an integer,
$a$, $b$ and $c$ be arbitrary parameters.
Then we have
\begin{align}
%
%
p_{n}(0;a,b,c,-c;q)
&=(-1)^ma^{m}b^mc^{2m}q^{m(3m-1)}(-c^2;q^2)_{m}
\nonumber\\&\times
p_m\left(x_{0};1,q,ab,-a^{-1}b^{-1}c^{-2}q^{-4m+2};q^2\right),
\label{eq:Psqrt-1even}\\
\intertext{if $n=2m$ is even, and}
%
%
p_{n}(0;a,b,c,-c;q)
&=(-1)^{m+1}a^{m}b^{m+1}c^{2m}(1+ab^{-1})q^{m(3m+1)}
(-c^2;q^2)_{m+1}
\nonumber\\&\times
p_m\left(x_{0};q,q^2,ab,-a^{-1}b^{-1}c^{-2}q^{-4m};q^2\right),
\label{eq:Psqrt-1odd}
\end{align}
if $n=2m+1$ is odd,
where $x_{0}=-\frac{ab^{-1}+a^{-1}b}2$.
\end{prop}
%
%
%
%
\begin{demo}{Proof}
We proceed by induction on $n$.
If $n=-1$ (resp. $n=0$), then 
the both sides of \eqref{eq:Psqrt-1odd} (resp. \eqref{eq:Psqrt-1even}) equals $0$ (resp. $1$).
Hence \eqref{eq:Psqrt-1even} and \eqref{eq:Psqrt-1odd} holds when $n=-1,0$.
%
Assume $n\geq1$ and \eqref{eq:Psqrt-1even} and \eqref{eq:Psqrt-1odd} hold when $n-1$ and $n-2$.
If $n=2m$, then 
\begin{align*}
&
p_{2m}(0;a,b,c,-c;q)
\intertext{by \eqref{eq:Askey-Wilson}}
&
=\frac{(ab,ac,-ac;q)_{2m}}{a^{2m}}\,{}_{4}\phi_{3}\left(
{{q^{-2m},-abc^2q^{2m-1},a\img,-a\img}\atop{ab,ac,-ac}};q,q
\right)
\intertext{by \eqref{eq:phi01}}
&
=\frac{(ab,ac,-ac;q)_{2m}}{a^{2m}}\,{}_{4}\phi_{3}\left(
{{q^{-2m+1},-abc^2q^{2m-2},a\img,-a\img}\atop{ab,ac,-ac}};q,q
\right)
\\&\qquad
-(1+a^2)(1+abc^2q^{4m-2})\frac{(abq,acq,-acq;q)_{2m-1}}{a^{2m}q^{2m-1}}
\\&\qquad\times
\,{}_{4}\phi_{3}\left(
{{q^{-2m+1},-abc^2q^{2m-1},aq\img,-aq\img}\atop{abq,acq,-acq}};q,q
\right)
\\&
=a^{-1}(1-abq^{2m-1})(1-acq^{2m-1})(1+acq^{2m-1})\,
p_{2m-1}(0;a,b,c,-c;q)
\\&\qquad
-a^{-1}(1+a^2)(1+abc^{2}q^{4m-2})\,
p_{2m-1}(0;aq,b,c,-c;q).
\intertext{By the induction hypothesis  \eqref{eq:Psqrt-1odd},
this equals}
&=(-1)^{m}a^{m-2}b^{m}c^{2m-2}q^{(m-1)(3m-2)}(-c^2;q^2)_{m}
\\&\times
\bigl\{
(1-abq^{2m-1})(1-acq^{2m-1})(1+acq^{2m-1})(1+ab^{-1})
\\&\qquad\qquad\times
p_{m-1}(x_{0};q,q^2,ab,-a^{-1}b^{-1}c^{-2}q^{-4m+4};q^2)
\\&\qquad
-q^{m-1}(1+a^2)(1+abc^{2}q^{4m-2})(1+ab^{-1}q)\,
\\&\qquad\qquad\times
p_{m-1}(x_{1};q,q^2,abq,-a^{-1}b^{-1}c^{-2}q^{-4m+3};q^2)\bigr\},
\end{align*}
where $x_{1}=-\frac{ab^{-1}q+a^{-1}bq^{-1}}2$.
Hence, by \eqref{eq:Askey-Wilson} again, this becomes
%
\begin{align*}
&=(-1)^{m}c^{2m}q^{m(3m-1)}(-c^2;q^2)_{m}(abq^2,abq^3,-c^{-2}q^{-4m+4};q^2)_{m-1}
\\&\times
\Biggl\{
-a^2(1-abq)(1+a^{-1}b)(1-a^{-2}c^{-2}q^{-4m+2})
\\&\qquad\times
\,
 {}_{4}\phi_{3}\left({
{q^{-2m+2},-c^{-2}q^{-2m+3},-a^2,-b^2}\atop{abq,abq^2,-c^{-2}q^{-4m+4}}
};q^2,q^2\right)
\\&\quad
-b^2(1+ab^{-1}q)(1+a^2)(1+a^{-1}b^{-1}c^{-2}q^{-4m+2})
\\&\qquad\times
\,
 {}_{4}\phi_{3}\left({
{q^{-2m+2},-c^{-2}q^{-2m+3},-a^2q^2,-b^2}\atop{abq^3,abq^2,-c^{-2}q^{-4m+4}}
};q^2,q^2\right)\Biggr\}.
\end{align*}
Using \eqref{eq:phi04b}, we obtain
\begin{align*}
&=(-1)^{m}c^{2m}q^{m(3m-1)}(-c^2,ab,abq,-c^{-2}q^{-4m+2};q^2)_{m}
\\&\qquad\times
\,
 {}_{4}\phi_{3}\left({
{q^{-2m},-c^{-2}q^{-2m+1},-a^2,-b^2}\atop{abq,ab,-c^{-2}q^{-4m+2}}
};q^2,q^2\right)
\intertext{by \eqref{eq:Askey-Wilson} again}
&=(-1)^{m}a^{m}b^{m}c^{2m}q^{m(3m-1)}
\left(-c^2;q^2\right)_{m}
\\&\qquad\times
\,
p_m\left(x_{0};q,1,ab,-a^{-1}b^{-1}c^{-2}q^{-4m+2};q^2\right).
\end{align*}
This proves \eqref{eq:Psqrt-1even} when $n=2m$.
%
%
%
%
%
%
%
It is well-known that the Askey-Wilson polynomials \eqref{eq:Askey-Wilson} are
symmetric with respect to the parameters $a$, $b$, $c$ and $d$.
Meanwhile, the expression of the Askey-Wilson polynomials by the ${}_{4}\phi_{3}$ series
apparently depends on the choice of the parameter $a$ in \eqref{eq:Askey-Wilson}.
Note that we have to choose this parameter carefully in each step
to apply the above contiguous relations.
Similarly, when $n=2m+1$, \eqref{eq:Psqrt-1odd} can be proven
using the definition \eqref{eq:Askey-Wilson} of the Askey-Wilson polynomials, 
the contiguous relations \eqref{eq:phi01}, \eqref{eq:phi02}
and the induction hypothesis \eqref{eq:Psqrt-1even} for $n=2m$.
We omit the details for the reader.
Hence we conclude that \eqref{eq:Psqrt-1even} and \eqref{eq:Psqrt-1odd} hold for arbitrary $m$.
\end{demo}
When $b=-a$, replacing $c$ by $b$, 
one gets incidentally the following known result due to Andrews (see \cite[(II.17)]{GR}).
%
%
\begin{corollary}
\begin{align*}
&
p_n(0;a,-a,b,-b;q)
\nonumber\\&=
\begin{cases}
(-1)^m(q,-a^2,-b^2,a^2b^2q^{2m};q^2)_m
&\mbox{ if $n=2m$,}\cr
0 &\mbox{ if $n=2m+1$. $\Box$}
\end{cases}
\end{align*}
\end{corollary}
%
%
%
%
\begin{demo}{Proof of Corollary~\ref{th:main}}
When $n=2m$ is even,
if we apply \eqref{eq:Psqrt-1even} to \eqref{eqtnabcn3aw} using
$
(b;q^2)_{m}\prod_{k=1}^{2m}(bq;q)_{k-2}=\prod_{k=1}^{m}\{(bq;q)_{2k-2}\}^2
$
then we obtain \eqref{eq:main-even2} by straightforward computation.
\eqref{eq:main-even} is derived from \eqref{eq:main-even2} by \eqref{eq:Askey-Wilson}
and using
$
(q;q^2)_m\prod_{k=1}^{2m}(q;q)_{k-1}=\prod_{k=1}^{m}\{(q;q)_{2k-1}\}^{2}
$,
$
(aq^{r+1};q^2)_{m}\prod_{k=1}^{2m}(aq;q)_{k+r+1}=\prod_{k=1}^{m}\{(aq;q)_{2k+r-1}\}^2
$
and
$
\frac{(a^{-1}b^{-1}q^{1-4m-r};q^2)_{m}}{\prod_{k=1}^{2m}(abq^2;q)_{k+2m+r-2}}
=\frac{(-1)^ma^{-m}b^{-m}q^{-m(3m+r)}}{\prod_{k=1}^{m}\{(abq^2;q)_{2(k+m)+r-3}\}^2}
$.
\par\smallskip
When $n=2m+1$ is even,
if we apply \eqref{eq:Psqrt-1odd} to \eqref{eqtnabcn3aw} using
\[
(b;q^2)_{m+1}\prod_{k=1}^{2m+1}(bq;q)_{k-2}
=\prod_{k=1}^{m+1}(bq;q)_{2k-2}\cdot\prod_{k=1}^{m}(bq;q)_{2k-2}
\]
then we obtain \eqref{eq:main-odd2}.
Finally, 
\eqref{eq:main-odd} is obtained from \eqref{eq:main-odd2} by using \eqref{eq:Askey-Wilson},
\begin{align*}
&(a^{-1}b^{-1}q^{-4m-r};q^2)_{m}
=(-1)^ma^{-m}b^{-m}q^{-m(3m+r+1)}(abq^{2m+r+2};q^2)_{m},
\\&
(q^3;q^2)_{m}\prod_{k=1}^{2m+1}(q;q)_{k-1}
=\frac1{1-q}\cdot\prod_{k=1}^{m+1}(q;q)_{2k-1}\prod_{k=1}^{m}(q;q)_{2k-1},
\\&
(aq^{r+2};q^2)_{m}\prod_{k=1}^{2m+1}(aq;q)_{k+r-1}
=\prod_{k=1}^{m+1}(aq;q)_{2k+r-2}\prod_{k=1}^{m}(aq;q)_{2k+r},
\\&
\frac{(abq^{2m+r+2};q^2)_{m}}
{\prod_{k=1}^{2m+1}(abq^2;q)_{k+2m+r-1}}
\\&\quad
=\frac1{\prod_{k=1}^{m+1}(abq^2;q)_{2k+2m+r-2}\prod_{k=1}^{m}(abq^2;q)_{2k+2m+r-2}}.
\end{align*}
This completes the proof of Corollary~\ref{th:main}.
\end{demo}
%
%
%

%% file: mehta04.tex
%
%
%
\section{A quadratic relation}
\label{sec:quadratic-relation}
In this section we use the Desnanot-Jacobi adjoint matrix theorem
to derive a quadratic relation between the Askey-Wilson polynomials
for a special values of parameters.
\par\smallskip
First we recall the reader a well-known theorem for matrix.
The following identity is known as the Desnanot-Jacobi adjoint matrix theorem \cite[Theorem~3.12]{Bre}
%
%
%
%
%
%
%
\begin{eqnarray}
\det A^{[2,n-1]}_{[2,n-1]}
\det A^{[n]}_{[n]}
=\det A^{[n-1]}_{[n-1]}
\det A^{[2,n]}_{[2,n]}
-\det A^{[n-1]}_{[2,n]}
\det A^{[2,n]}_{[n-1]}.
\label{eq:Desnanot-Jacobi2}
\end{eqnarray}
Let
\[
D_{n}(a,b,c;q)
=\det\left(
(q^{i-1}-cq^{j-1})\frac{(aq;q)_{i+j-2}}{(abq^2;q)_{i+j-2}}
\right)_{1\leq i,j\leq n}
\]
and apply \eqref{eq:Desnanot-Jacobi2} to this determinant.
Then we obtain
\begin{align}
&D_{n}(a,b,c;q)D_{n-2}(aq^2,b,c;q)
=\frac{q(aq;q)_{2}}{(abq^2;q)_{2}}
\cdot
D_{n-1}(a,b,c;q)D_{n-1}(aq^2,b,c;q)
\nonumber\\&\quad
-\frac{q(1-aq)^{n}(1-abq^3)^{n-2}}{(1-aq^2)^{n-2}(1-abq^2)^{n}}
\cdot
D_{n-1}(aq,b,cq;q)D_{n-1}(aq,b,cq^{-1};q).
\label{eq:DJ}
\end{align}
Hence we can substitute \eqref{eqtnabcn3aw} into \eqref{eq:DJ}.
Then we obtain
\begin{align}
&
aq(1-q^{n-1})(1-bq^{n-2})\cdot
p_{n}(0;a^{\frac12}c^{\frac12}q^{\frac12}\img,-a^{\frac12}c^{-\frac12}q^{\frac12}\img,b^{\frac12}\img,-b^{\frac12}\img;q)
\nonumber\\&\qquad\times
p_{n-2}(0;a^{\frac12}c^{\frac12}q^{\frac32}\img,-a^{\frac12}c^{-\frac12}q^{\frac32}\img,b^{\frac12}\img,-b^{\frac12}\img;q)
\nonumber\\&=
(1-aq^{n})(1-abq^{n})\cdot
p_{n-1}(0;a^{\frac12}c^{\frac12}q^{\frac12}\img,-a^{\frac12}c^{-\frac12}q^{\frac12}\img,b^{\frac12}\img,-b^{\frac12}\img;q)
\nonumber\\&\qquad\times
p_{n-1}(0;a^{\frac12}c^{\frac12}q^{\frac32}\img,-a^{\frac12}c^{-\frac12}q^{\frac32}\img,b^{\frac12}\img,-b^{\frac12}\img;q)
\nonumber\\&-
(1-aq)(1-abq^{2n-1})\cdot
p_{n-1}(0;a^{\frac12}c^{\frac12}q^{\frac32}\img,-a^{\frac12}c^{-\frac12}q^{\frac12}\img,b^{\frac12}\img,-b^{\frac12}\img;q)
\nonumber\\&\qquad\times
p_{n-1}(0;a^{\frac12}c^{\frac12}q^{\frac12}\img,-a^{\frac12}c^{-\frac12}q^{\frac32}\img,b^{\frac12}\img,-b^{\frac12}\img;q).
\label{eq:quadratic-relation0}
\end{align}
Replacing 
$a^{\frac{1}{2}}c^{\frac{1}{2}}q^{\frac{1}{2}}\img$, 
$-a^{\frac{1}{2}}c^{-\frac{1}{2}}q^{\frac{1}{2}}\img$ 
and
$b^{\frac{1}{2}}\img$ 
by $a$, $b$ and $c$, respectively
in \eqref{eq:quadratic-relation0},
we obtain the following corollary.
%
%
%
%
%
%
\begin{corollary}
\label{cor:quadratic-relation}
Let $n$ be a positive integer and $a$, $b$, $c$ and $q$ parameters.
Then we have
\begin{align}
&ab(1-q^{n-1})(1+c^2q^{n-2})
p_n(0;a,b,c,-c;q)p_{n-2}(0;aq,bq,c,-c;q)\nonumber\\
&=
(1-abq^{n-1})(1+abc^2q^{n-1})
p_{n-1}(0;a,b,c,-c;q)p_{n-1}(0;aq,bq,c,-c;q)\nonumber\\
&\quad
-
(1-ab)(1+abc^2q^{2n-2})
p_{n-1}(0;aq,b,c,-c;q)p_{n-1}(0;a,bq,c,-c;q).
\label{eq:quadratic-relation}
\end{align}
In other word
\begin{align}
&abq(1-q^{n-1})(1+c^2q^{n-2})
{}_4\phi_3\left(
\begin{matrix}
q^{-n},-abc^2q^{n-1},a\img,-a\img\cr
ab,ac,-ac
\end{matrix}
;q,q\right)
\nonumber\\&\qquad\times
{}_4\phi_3\left(
\begin{matrix}
q^{-n+2},-abc^2q^{n-1},aq\img,-aq\img\cr
abq^2,acq,-acq
\end{matrix}
;q,q\right)\nonumber\\
&=
(1-abq^n)(1+abc^2q^{n-1})
{}_4\phi_3\left(
\begin{matrix}
q^{-n+1},-abc^2q^{n-2},a\img,-a\img\cr
ab,ac,-ac
\end{matrix}
;q,q\right)
\nonumber\\&\qquad\times
{}_4\phi_3\left(
\begin{matrix}
q^{-n+1},-abc^2q^n,aq\img,-aq\img\cr
abq^2,acq,-acq
\end{matrix}
;q,q\right)\nonumber\\
&\quad
-(1-abq)(1+abc^2q^{2n-2})
{}_4\phi_3\left(
\begin{matrix}
q^{-n+1},-abc^2q^{n-1},aq\img,-aq\img\cr
abq,acq,-acq
\end{matrix}
;q,q\right)
\nonumber\\&\qquad\times
{}_4\phi_3\left(
\begin{matrix}
q^{-n+1},-abc^2q^{n-1},a\img,-a\img\cr
abq,ac,-ac
\end{matrix}
;q,q\right).
\label{eq:quadratic-relation2}
\end{align}
\end{corollary}
%
%
%
%
%
%
%
%
%
%
%
%
Here we derive Corollary~\ref{cor:quadratic-relation}
as a corollary of Theorem~\ref{th:corthm2}.
Note that, if one can prove the quadratic relation
\eqref{eq:quadratic-relation} or \eqref{eq:quadratic-relation2} directly,
then he can prove Theorem~\ref{th:corthm2} using the Desnanot-Jacobi adjoint matrix theorem
\eqref{eq:Desnanot-Jacobi2}.
This was our first strategy to prove Theorem~\ref{th:corthm2},
but we have found it is not so easy to prove the quadratic relation.
Hence we prove Theorem~\ref{lem:thm1} which is more general,
and then derive Theorem~\ref{th:corthm2}.
\par\smallskip
Finally,
 we state a conjecture which generalizes \eqref{eq:quadratic-relation}.
We can observe that the following quadratic equation with two more parameters could hold concerning to the Askey-Wilson polynomials:
%
%
%
%
%
%
\begin{conjecture}
Let $n$ be a positive integer and $x$, $a$, $b$, $c$, $d$ and $q$ parameters.
Then we have
\begin{align}
&ab(1-q^{n-1})(1-cdq^{n-2})
p_n(x;a,b,c,d;q)p_{n-2}(x;aq,bq,c,d;q)\nonumber\\
&=
(1-abq^{n-1})(1-abcdq^{n-1})
p_{n-1}(x;a,b,c,d;q)p_{n-1}(x;aq,bq,c,d;q)\nonumber\\
&\quad
-
(1-ab)(1-abcdq^{2n-2})
p_{n-1}(x;aq,b,c,d;q)p_{n-1}(x;a,bq,c,d;q).\label{eqmw3}
\end{align}
\end{conjecture}
%
%
%
%
\par\medskip\noindent
{\large\bf Concluding Remarks}
This conjecture may hint us there could exist a more general formula
than Theorem~\ref{th:corthm2}.
But it is not an easy task to find the appropriate entry of the determinant
which gives this quadratic relation.
\par
It is also an interesting problem to find 
a combinatorial application of Theorem~\ref{th:corthm2}.
If $c=0$, then \eqref{eq:q-general-hankel} is the generating function
of the Dyck paths with certain weights (see \cite{ITZ1}).
For $c=1$, the Pfaffian \eqref{eq:pf-special} can enumerate certain reverse plane partitions (see \cite{ITZ2}).
We have been trying to find an application of Theorem~\ref{th:corthm2}
and it is not so easy to find an appropriate lattice path and its weights.
It will be left for the future work.
%
%
%
%
%
%
%